%% file: NAHS_PII_ver2.tex
\documentclass[onefignum,onetabnum]{siamart171218}

\usepackage{amssymb}
\usepackage{amsopn}

\usepackage{tikz}
\usetikzlibrary{scopes}
\usetikzlibrary{arrows}
\usetikzlibrary{arrows.meta}
\usetikzlibrary{decorations.pathreplacing,decorations.markings}
\usepackage{lipsum}
\usepackage{amsfonts}
\usepackage{graphicx}
\usepackage{epstopdf}
\usepackage{algorithmic}
\usepackage{centernot}
\usepackage{color}
\ifpdf
\DeclareGraphicsExtensions{.eps,.pdf,.png,.jpg}
\else
\DeclareGraphicsExtensions{.eps}
\fi

\title{Algorithms for optimal control of hybrid systems with
sliding modes}

\author{Rados\l aw Pytlak\thanks{Faculty of Mathematics and Information Science, Warsaw University of Technology, 00-662 Warsaw, Poland 
  (\email{radoslaw.pytlak@pw.edu.pl}.}
\and Damian Suski\thanks{Institute of Automatic Control and Robotics, Warsaw University of Technology, 02-525 Warsaw, Poland 
  (\email{damian.suski@pw.edu.pl}.}}

\begin{document}

\maketitle

\begin{abstract}
This paper concerns two algorithms for solving optimal control problems with hybrid systems. The first algorithm aims at hybrid systems exhibiting sliding modes. The first algorithm has several features which distinguishes it from the other algorithms for problems described by hybrid systems. First of all, it can cope with hybrid systems which exhibit sliding modes. Secondly, the systems motion on the switching surface is described by index 2 differential--algebraic equations and that guarantees accurate tracking of the sliding motion surface. Thirdly, the gradients of the problems functionals are evaluated with the help of adjoint equations. The adjoint equations presented in the paper take into account sliding motion and exhibit jump conditions at transition times. We state optimality conditions in the form of the weak maximum principle for optimal control problems with hybrid systems exhibiting sliding modes and with piecewise differentiable controls. The second algorithm is for optimal control problems with hybrid systems which do not exhibit sliding motion. In the case of this algorithm we assume that control functions are ${\mathcal L}^\infty$ measurable functions. For each algorithm, we show that every accumulation point of the sequence generated by the algorithm satisfies the weak maximum principle.
\end{abstract}



\input{sec_introductionSC_P2.tex}
\input{sec_algorithm_exactPenalty_SCP2.tex}
\input{sec_adjointEqsSC_P2.tex}

\input{sec_withoutSMSC_P2.tex}

\input{sec_conclusionsSC_P2.tex}

\bibliographystyle{siamplain}
\bibliography{NAHS_p2}
\input{sec_appendixSC_P2_ver2.tex}

\end{document}

%% file: sec_introductionSC_P2.tex
\section{Introduction}

Hybrid systems are systems with mixed discrete-continuous dynamics (\cite{ss2000},\cite{bbm1999}). The set of discrete states $ \mathcal{Q} $ consists of finite number of elements denoted by $ q $. The admissible controls set $ \mathcal{U} $ consists of measurable functions $ u:I \rightarrow U $  defined on a closed interval $ I $ with the values in $ U \in \mathbb{R}^m $. The continuous dynamics in each discrete state is described by ordinary differential equations (ODEs)
\begin{eqnarray}
{\displaystyle x' = f(x,u)}\label{explicit}
\end{eqnarray}
or more generally by differential--algebraic equations (DAEs)
\begin{eqnarray}
&{\displaystyle 0 = F(x',x,u),}\label{implicit}
\end{eqnarray}
where $ x \in \mathbb{R}^n $, $ f: \mathbb{R}^n \times U \rightarrow \mathbb{R}^n $, $ F: \mathbb{R}^n \times \mathbb{R}^n \times U \rightarrow \mathbb{R}^n $.
The transitions between discrete states are triggered when the condition of the form 
\begin{equation}
\eta(x,u) \leq 0
\end{equation}
stops to be satisfied, where  $ \eta: \mathbb{R}^n \times U \rightarrow \mathbb{R} $. The functions $ \eta(x,u) $ are called {\it guards}. We restrict our analysis to systems with autonomous transitions and without state jumps during transitions. It should be noted that we consider only the sliding motion surfaces of co--dimension 1, i.e. we do not consider the case, when the sliding motion surfaces may intersect (\cite{gh2021},\cite{kk2019}).

The optimal control problem with a hybrid system has been considered in many papers. The necessary optimality conditions for a class of hybrid systems without state jumps have been first formulated in \cite{w1966}. In \cite{bh1975} the variational methods have been used to formulate adjoint equations for systems with state jumps. The Pontryagin maximum principle for hybrid systems with state jumps has been formulated for several classes of hybrid systems in \cite{s1999}, \cite{s2004}, \cite{sc2007}, \cite{sc2009}, \cite{tc2010}, \cite{tc2011}, \cite{tc2013}, \cite{pc2013}. In papers  \cite{s2004}, \cite{sc2007}, \cite{sc2009}, \cite{tc2010}, \cite{tc2013} also algorithms based on the hybrid maximum principle are discussed. Our paper does not consider optimal control problems in which switching times are decision variables (\cite{scgb06}).

In none of these papers an optimal control problem with hybrid systems exhibiting sliding modes has been considered. In \cite{ps2018} we have introduced such an algorithm based on the adjoint equations for differential--algebraic equations describing dynamics of hybrid system with sliding motion. \cite{ps2018} provides also a computational example which shows the validity of the approach. However, \cite{ps2018} does not contain the derivation of the adjoint equations and a proof of global convergence of the proposed optimization algorithm. The reason for that is that in order to do that  it is necessary to provide sensitivity analysis of trajectories of hybrid systems exhibiting sliding modes. That analysis was stated in \cite{ps2020a} and is used in this paper to derive adjoint equations and a weak version of the maximum principle provided that controls are represented by piecewise smooth controls and the proposed constraint qualification holds. Furthermore, we propose an algorithm which generates a sequence of approximating optimal controls which converges to a control which satisfies optimality conditions stated by the weak maximum principle.  

The algorithm, convergence analysis and the derivation of the necessary optimality conditions are based on an exact penalty function associated with objective function and functions defining constraints of the considered problem. The paper emphasizes  the role of the proposed constraint qualification in the process of solving  optimal control problem with hybrid systems.

The main focus of the paper is on the problems with hybrid systems which can exhibit sliding motion. In that case we assume that controls are piecewise smooth functions. In the paper we also discuss the method which is intended for control problems with hybrid systems whose trajectories do not experience sliding motion. In the case of these optimal control problems we assume that control functions are elements of the space ${\mathcal L}^\infty$.

The convergence is established under the assumption that continuous time system equations are given. In the accompanying papers (\cite{ps19a},\cite{ps19b}), which are extension of the paper \cite{ps2021}, we show that the algorithm convergence can still be attained if state trajectories are obtained by an implicit Runge--Kutta method provided that the integration procedure step sizes go to zero. One of the applications of the method was the procedure for planning a haemodialysis process (\cite{PytlakSuskiCDC}) which was one of the motivations for carrying out the research presented in this paper.

The computational approach to hybrid systems we advocate in this paper (and in \cite{ps17}, \cite{ps19a}, \cite{ps19b}) assumes that for a given control function we attempt to follow true systems trajectories as close as possible by using integration procedures with a high order of convergence. It is in contrast to the approach which uses approximations to the discontinuous right-hand side for the differential equations or inclusions by a smooth right--hand side (\cite{sa2010}, \cite{sa2012}, \cite{mhrl2015}). The smoothing approach demonstrates that it can be used to find approximate solutions to optimal control problems with systems which exhibit sliding modes however there is still a lack of evidence that it is capable of efficiently finding optimal controls and state trajectories with high accuracy by choosing proper values of smoothing parameters.

%% file: sec_algorithm_exactPenalty_SCP2.tex
\section{Calculating the optimal control}
\label{CalcOptCon}

In our work we consider hybrid systems with sliding modes. For the sake of simplicity let us consider a hybrid system with two discrete states collected in a set $ \mathcal{Q} = \{1,2 \}$. Let us assume that the invariant sets are  $ \mathcal{I}(1) = \left\{ x\in \mathbb{R}^n: h(x)\leq 0 \right\} $ and $ \mathcal{I}(2) = \left\{ x\in \mathbb{R}^n: h(x)\geq 0 \right\} $ where $ h: \mathbb{R}^n \rightarrow \mathbb{R} $. If the hybrid system starts its evolution from a discrete state $ q = 1 $ the continuous state evolves according to an equation $ x' = f_1(x,u) $. At a transition time $ t_t $ the continuous state trajectory reaches the boundary of an invariant set and we have $ h(x(t_t)) =0$. The first order condition which guarantees that the continuous state trajectory will leave the invariant set $ \mathcal{I}(1) $ is (\cite{ljs2000})
\begin{equation}
h_x(x(t_t))f_1(x(t_t),u(t_t)) > 0 \label{sw1}
\end{equation}
where $ h_x^T(x) $ is the normal vector to a surface 
\begin{equation}
\Sigma = \{ x\in \mathbb{R}^n: h(x) = 0 \}
\end{equation} at $ x $. If at a transition time we also have 
\begin{equation}
h_x(x(t_t))f_2(x(t_t),u(t_t)) > 0 \label{sw2}
\end{equation}
then the discrete state changes from  $ q=1 $ to $ q = 2 $ and the continuous state continues the evolution according to the equation $ x' = f_2(x,u) $. If at a transition time we have 
\begin{equation}
h_x(x(t_t))f_2(x(t_t),u(t_t)) < 0
\end{equation}
then both vector fields $ f_1(x,u) $ and $ f_2(x,u) $ point towards the surface $ \Sigma $ and we face the {\it sliding motion} phenomenon (\cite{dl2009}).

The sliding motion can be handled with the concept of {\it Filippov solutions} (\cite{fil88},\cite{dl2009}). We say that a continuous state trajectory is a Filippov solution of the considered hybrid system if
\begin{eqnarray}
&{\displaystyle 
	x' = \left \{
	\begin{array} {ll}
	f_1(x,u) & {\rm if}\ h(x) < 0 \\
	\overline{\rm co} \left ( f_1(x,u), f_2(x,u) \right ) & {\rm if}\ h(x) = 0 \\
	f_2(x,u) & {\rm if}\ h(x) > 0 
	\end{array} \right .}\nonumber 
\end{eqnarray}
for almost all t from its domain of definition. Here, $ \overline{\rm co} \left ( f_1(x,u), f_2(x,u) \right ) $ is is the minimal closed convex set containing $f_1$ and $f_2$, that is
\begin{eqnarray}
&{\displaystyle 
\overline{\rm co} \left ( f_1, f_2 \right )= \left \{ f_F\in {\mathcal R}^n:\  f_F = f_1 + \alpha \left (f_2-f_1\right ),\ \alpha \in [0,1] \right \}.}\label{fFDef}
\end{eqnarray}

During the sliding motion the continuous state trajectory must stay in the surface $ \Sigma $, so the condition
\begin{equation}
h_x(x)f_F(x,u)=0 \label{fFTangency}
\end{equation} 
must be satisfied. From (\ref{fFDef}) and (\ref{fFTangency}) it is easy to find the formula for $ \alpha $ coefficient
\begin{equation}
\alpha(x,u) = \frac{h_x(x)f_1(x,u)}{h_x(x)\left(f_1(x,u)-f_2(x,u) \right)}.\label{AlphaF}
\end{equation} 

When dealing with hybrid systems we have to pay special attention to systems behavior at switching times $t_t$, on the left and the right of these points.
To this end we define
\begin{eqnarray}
&{\displaystyle 
	u(t_t^-) =  \lim_{t\rightarrow t_t, t<t_t} u(t),\  
	u(t_t^+) =  \lim_{t\rightarrow t_t, t>t_t} u(t),} \label{Jump2}
\end{eqnarray}
and we assume that $u$ is continuous from the left which means that 
\begin{eqnarray}
&{\displaystyle 
u(t_t^-) = u(t_t)\ \ {\rm and}\ u(t_t^+) \neq u(t_t)} \label{Jump4}
\end{eqnarray}
in general. 

This notation at a neighborhood of time $t_t$ apply also to other functions, for example
\begin{eqnarray}
&{\displaystyle 
	x(t_t^-) =  \lim_{t\rightarrow t_t, t <t_t} x(t),\ 
	x(t_t^+) = \lim_{t\rightarrow t_t, t >t_t} x(t).} \label{Jump2x}
\end{eqnarray}

In this setting general conditions (\ref{sw1})-(\ref{sw2}) for changing discrete states must be stated as follows:
\begin{eqnarray}
h_x(x(t_t^{-}))f_1(x(t_t^{-}),u(t_t^{-})) & > & 0 \label{sw1a} \\
h_x(x(t_t^+))f_2(x(t_t^+),u(t_t^+)) & > & 0. \label{sw2a}
\end{eqnarray}

If the system is in the sliding mode and at the time $t_t$ transits, for example, to state $q=1$, the conditions for changing discrete state can be stated as follows (then as a guard $\eta$ we take $\alpha$)
\begin{eqnarray}
\alpha(x(t_t^-),u(t_t^-)) &=& 0 \label{sw3lsc} \\
h_x(x(t_t^+))f_1(x(t_t^+),u(t_t^+)) & < & 0 \label{sw3alsc} \\
h_x(x(t_t^+))f_2(x(t_t^+),u(t_t^+)) & < & 0 \label{sw3blsc}
\end{eqnarray}

The introduced convention is also applied to  integrals, for example
\begin{eqnarray}
&{\displaystyle \int_{t_0}^{t_t^-}p(x(t),u(t))dt = \lim_{t\rightarrow t_t,t< t_t}\int_{t_0}^t p(x(\tau),u(\tau))d\tau. }\nonumber
\end{eqnarray}

But taking piecewise smooth functions as controls requires also checking at each point $t_j$ (see (\ref{pwc2})) whether the system is still at a sliding mode after control changes at time $t_j$, if not a new discrete state at this time must be determined. For lack of space we do take into account these nonautonomous switches times. However, our analysis could be easily extended with this respect under additional assumption ${\bf (H2)}\ iv)$ stated in \cite{ps2020a}. 

During the sliding motion the continuous state trajectory should stay in a set $ \Sigma $. To keep this condition satisfied during the numerical integration, we follow the approach proposed in \cite{ap1998} and integrate the differential--algebraic equations (DAEs)
\begin{eqnarray}
x' &=& f_F(x,u) + h_x(x)^Tz \label{slidingDAEDiff} \\
0 &=& h(x) \label{slidingDAEAlg}
\end{eqnarray}
instead of ODEs
\begin{equation}
x' = f_F(x,u). \label{Fil2}
\end{equation}

It can be shown that at a sliding motion equations (\ref{slidingDAEDiff})--(\ref{slidingDAEAlg}) and (\ref{Fil2}) are equivalent (\cite{ps2020a}). We prefer using (\ref{slidingDAEDiff})--(\ref{slidingDAEAlg}) instead of (\ref{Fil2}) for the reason of getting more stable numerical solution on the switching surface.

The sliding motion can be treated as the third discrete state with $ q = 3 $. Treating the sliding motion as another discrete state implies that guards functions $h$ should depend also on control variables $u$ and system equations should be given, in general, by implicit equations (\ref{slidingDAEDiff})-(\ref{slidingDAEAlg}).

From now on  the invariant set at a discrete state $q$ will be defined as follows
\begin{eqnarray}
&{\displaystyle 
\mathcal{I}(q)  = \left\{ (x,z,u)\in \mathbb{R}^{n+1+m}: \eta_q(x,u)\leq 0 \right\},\ \eta_q: \mathbb{R}^{n} \times U\rightarrow \mathbb{R}^{n_{{\mathcal I}_q}}.}
\label{defInvariantSeta}
\end{eqnarray}
and $n_{{\mathcal I}_q}$ denotes the number of functions defining the invariant set for state $q$.

Furthermore, the system evolution is related to the functions $f_q$ and $h_q$ defining, for a given discrete state $q$, DAEs:
\begin{eqnarray}
x' & = & f_q (x,u) + \left (h_q\right )_x(x)^Tz \nonumber \\
0  & = & h_q(x).\nonumber 
\end{eqnarray}
In particular, we can have $h_q (x)\equiv 0$ for some $q$ and then the guarding function $\eta_q$ is not dependent on $u$. 

If we assume that a switching surface function $\eta$ depends on $u$ we need continuity from the left of the function $\eta(x(t),u(t))$ to locate uniquely  a switching time $t_t$. For that reason we also assume that control functions $u$ are piecewise smooth.


Consider the system
\begin{eqnarray}
x' & = & f(x,u) + h_x(x)^Tz = f_r(x,z,u) \label{NewslidingDAEDiffInd1} \\
0  & = & h(x)\label{NewslidingDAEAlgInd1}
\end{eqnarray}
on $[t_0,t_f]$. We assume that on some subintervals $[t^0,t^1]$ we can have $h(x)\equiv 0$ and in that case system (\ref{NewslidingDAEDiffInd1})--(\ref{NewslidingDAEAlgInd1}) reduces to ODEs. Since we consider, for the simplicity of presentation, that our hybrid system can be only in three states $q=1,2,3$ (and the third state reserved for a sliding motion, in that state $f\equiv f_F$) our state trajectory on the entire horizon $[t_0,t_f]$ will consist of trajectories of type $(x,z)$ on subintervals $A^s_i$, $i\in I_s$, and of type $x$ on subintervals $A^1_i$, $i\in I_1$ (if trajectory is defined by the function $f_1$) and subintervals $A^2_i$, $i\in I_2$ (if trajectory is defined by the function $f_2$). We have $\cup_{i\in I^s_i}A^s_i\cup_{i\in I_1}A^1_i\cup_{i\in I_2}A^2_i = [t_0,t_f]$. 

Our approach is based on sensitivity analysis stated in \cite{ps2020a} where it is shown that solutions to linearized equations of DAEs have properties similar to those exhibited by linearized equations to ordinary differential equations, provided that some conditions are met. Namely, if we denote by $(x,z)$ the solution of (\ref{NewslidingDAEDiffInd1})--(\ref{NewslidingDAEAlgInd1}) for given $u$ (in some situations, if it does not lead to the confusion, we will also write $(x^u,z^u)$ for $(x,z)$), by $(x^{d},z^{d})$ the solution for the control $u+d$, and by $(y^{x,d},y^{z,d})$ the solution to the linearized equations:
\begin{eqnarray}
\left (y^{x}\right )' &=& f_x(x,u)y^x + \left (h_x(x)^Tz\right)_x y^x + h_x(x)^T y^z + f_u(x,u) d \label{LslidingDAEDiff} \\
0 &=& h_x(x) y^x, \label{LslidingDAEAlg}
\end{eqnarray}
then the following will be satisfied
\begin{eqnarray}
\|(x,z)\|_{{\mathcal L}^\infty} & \leq & c_1 \label{Lin1} \\
\| (x^{d},z^{d}) - (x,z)\|_{{\mathcal L}^\infty} & \leq & c_2 \|d\|_{{\mathcal L}^\infty} \label{Lin2} \\
\| (y^{x,d},y^{z,d})\|_{{\mathcal L}^\infty} & \leq & c_3 \|d\|_{{\mathcal L}^1} \label{Lin3}
\end{eqnarray}
for some positive constants $c_1$, $c_2$, $c_3$. Furthermore, for any $\varepsilon$ satisfying $\varepsilon_t > \varepsilon > 0$, where $\varepsilon_t$ is as specified in ${\bf (H3)}$, there exists function $o:[0,\infty]\rightarrow (0,\infty)$ such that $\lim_{s\rightarrow 0^{+}} o(s)/s = 0$ and
\begin{eqnarray}
&{\displaystyle 
	\| x^{d}(t) -  \left (x(t) + y^{x,d}(t)\right )\|  \leq  o\left ( \|d\|_{{\mathcal L}^\infty}\right ),}\nonumber \\
&{\displaystyle 
\forall t\in [t_0,t_t^1-\varepsilon] \bigcup [t_f-\varepsilon, t_f]  
	\bigcup_{i\in I_t(u)\backslash \{N_t(u)\}} [t_t^i,t_t^{i+1}-\varepsilon].}\label{Lin4}
\end{eqnarray}
Here, the points $\{t_t^i\}_{i=1}^{N_t(u)}$ (and the related to them the set of indices $I_t(u)$) are as specified in the assumption ${\bf (H3)}$.

In order to specify the relations for differentials of $x$ at a switching point $t_t$, denoted by $dx^d$, additional notation has to be introduced. Suppose that a control $u$ is perturbed by $d$, then a new switching point $t_t^d$ will occur. In general setting, we assume that before  switching of discrete states the system evolves according to equations (\ref{NewslidingDAEDiffInd1}--(\ref{NewslidingDAEAlgInd1}) with $f$ and $h$ replaced by $f_1$ and $h_1$ respectively and after the switching the system is described by (\ref{NewslidingDAEDiffInd1}--(\ref{NewslidingDAEAlgInd1}) with $f_2$ and $h_2$ instead of $f$ and $h$. Then, we have
\begin{eqnarray}
\left \| x^d_1(t_t^{d-})-x_1(t_t^-)- d\left (x_1\right)^d_-\right \| & \leq & o(\|d\|_{{\mathcal L}^\infty}), \label{Lin5} \\
\left \| x^d_2(t_t^{d+})-x_2(t_t^+)- d\left (x_2\right)^d_+\right \| & \leq & o(\|d\|_{{\mathcal L}^\infty}), 
\label{Lin6}
\end{eqnarray}
for all $u$ and $d$ such that $u+d\in {\mathcal U}$ where $d\left (x_1\right)^d_- = y^{x,d}_1(t_t^-) + f_{1r}(x_1(t_t^-),z_1(t_t^-),$ $u(t_t^-))dt_t^d$, $x_1$, $x^d_1$, $z_1$ are solutions to equations (\ref{NewslidingDAEDiffInd1})--(\ref{NewslidingDAEAlgInd1}) with $f=f_1$, $h=h_1$, $y^{x,d}_1$ are solutions to the corresponding linearized equations associated with equations (\ref{NewslidingDAEDiffInd1})--(\ref{NewslidingDAEAlgInd1}),
$d\left (x_2\right)^d_+ = y^{x,d}_2(t_t^+) + f_{2r}(x_2(t_t^+),z_2(t_t^+),u(t_t^+))dt_t^d$,
$x_2$, $x^d_2$, $z_2$ are solutions to equations (\ref{NewslidingDAEDiffInd1})--(\ref{NewslidingDAEAlgInd1}) with $f=f_2$, $h=h_2$ and $y^{x,d}_2$ are solutions to the corresponding linearized equations associated with equations (\ref{NewslidingDAEDiffInd1})--(\ref{NewslidingDAEAlgInd1}). Here, $dt_t^d$ is the differential of the switching point $t_t$ and $o$ are such that $\lim_{s\rightarrow 0^+}|o(s)|/s =0$; $dx^d_-$ and $dx^d_+$ differentials correspond to systems behavior before and after the switching ($dx^d_-=dx^d_+$).

In order to present the conjectures under which (\ref{Lin1})--(\ref{Lin4}) hold we introduce the function $e(x,z,u) = h_x(x)f(x,u) + \|h_x(x)\|^2z$. Then under {\bf (H1)}, {\bf (H2)} and {\bf (H3)} the relations (\ref{Lin1})--(\ref{Lin4}) are satisfied.
\begin{description}
	\item[${\bf (H1)}$]  $U$ is a compact convex set in ${\mathbb R}^m$, $f: {\mathbb R}^n\times {\mathbb R}^m\rightarrow {\mathbb R}^n$, $h: {\mathbb R}^n\rightarrow {\mathbb R}$, \\
	i) function $f(\cdot,\cdot)$ is differentiable, $f_x$ $f_u$ are continuous and there exists  $0<k_f^1 < +\infty$ such that 
	\begin{eqnarray}
	&{\displaystyle 
		\| f_x(x,u) \| \leq k_f^1} \label{SH2b}
	\end{eqnarray}
	for all $(x,u)\in {\mathbb R}^n\times U$,\\
	ii) if $h(x)\not\equiv 0$ then $(x(t_0),z(t_0))$ are such that
	\begin{eqnarray}
	&{\displaystyle   h(x(t_0)) =0,\ z(t_0)=0,\  h_x(x(t_0))f(x(t_0),u(t_0))=0,} \label{initial1}
	\end{eqnarray}
	$h(\cdot)$ is twice differentiable and $\nabla^2 h$  is continuous. Moreover, there exist $\varepsilon > 0$, $ 0 < k_{hf} <+\infty$, $0< k_e^1 < +\infty$, $0 < k_f^2 < +\infty$, $0< k_f^3 < +\infty$, $0 < k_h^1 < k_h^2 < +\infty$, $0 < k_h^3 < +\infty$ with the properties
	\begin{eqnarray}
	&{\displaystyle \| f(x,u) \| \leq k_f^2} \label{SM2ba}\\
	&{\displaystyle \| f_u(x,u) \| \leq k_f^3} \label{SM2bc} \\
	&{\displaystyle 
		\|h_x(x)f(x,u)\| \leq k_{hf}} \label{SH2a} \\
	&{\displaystyle 
		\| e_{x,z,u}(\hat{x},\hat{z},\hat{u}) - e_{x,z,u}(x,z,u)\|\leq k_e^1\|(\hat{x},\hat{z},\hat{u}) - (x,z,u) \|} \label{SH2c} \\
	&{\displaystyle 
		k_h^1 \leq \|h_x(x) \| \leq k_h^2} \label{SH2e} \\
	&{\displaystyle 
		\|\nabla^2 h(x)\| \leq k_h^3} \label{SH2f}
	\end{eqnarray}
	for all  $(x,z,u)$, $(\hat{x},\hat{z},\hat{u})$ in $\mathbb{R}^n\times \mathbb{R} \times U$ such that $|h(x)|\leq\varepsilon$, $| h(\hat{x})|\leq \varepsilon$.
\end{description}

In ${\bf (H1)}$ it is assumed that $U$ is a compact convex set. In order to show (\ref{Lin1})--(\ref{Lin4}) it is sufficient to assume that $U$ is bounded, but the convergence analysis which is based on (\ref{Lin1})--(\ref{Lin4}) requires the stronger assumption.

\begin{description}
	\item[${\bf (H2)}$]  i) Function $\eta(\cdot,\cdot)$ is differentiable and there exist $0<L_1<+\infty$ and $0<L_2 <+\infty$ such that
	\begin{eqnarray}
	\left \| \eta_{x,u}(\hat{x},\hat{u}) - \eta_{x,u}(x,u)\right \| & \leq  & L_1\|(\hat{x},\hat{u}) - (x,u) \| \label{H3a} \\
	\left | f_{r}(\hat{x},\hat{z},\hat{u}) -  f_{r}(x,z,u)\right | & \leq & L_2\|(\hat{x},\hat{z},\hat{u}) - (x,z,u) \|      \label{H3b}	
	\end{eqnarray}
	for all  $(x,z,u)$, $(\hat{x},\hat{z},\hat{u})$ in $\mathbb{R}^n\times \mathbb{R} \times U$.
	\vspace{2mm}
	
	ii) There exists $\varepsilon > 0$ and $0<L_3 <+\infty$ such that for all switching points $t_t$ (each switching time corresponds to some $u\in {\mathcal U}$) and for all their perturbations $t_t^d$ triggered by perturbations $d$ such that $u+d\in {\mathcal U}$, $\|d\|_{{\mathcal L}^\infty}\leq \varepsilon$, and for all $\theta\in [0,1]$ we have
	\begin{eqnarray}
	&{\displaystyle 
		\left |\eta_x(x(\tau^{\theta ,d}(t_t)),u(\tau^{\theta ,d}(t_t) )) f_{r}(x(\tau^{\theta ,d}(t_t)),z(\tau^{\theta ,d}(t_t)),u(\tau^{\theta ,d}(t_t) ))\right . }\nonumber \\ &{\displaystyle \left . + \eta_u(x(\tau^{\theta ,d}(t_t)),u(\tau^{\theta ,d}(t_t) ))u'(\tau^{\theta ,d}(t_t))\right |\geq L_3} 
	\label{H3c}	
	\end{eqnarray}
	where
	\begin{eqnarray}
	\tau^{\theta ,d}(t_t) & = & t_t + \theta (t^d_t - t_t). 
	\label{Midt} 
	\end{eqnarray}
	\vspace{2mm}
	iii) Perturbations $d\in {\mathcal U}-u$ are such that 
	\begin{eqnarray}
	&{\displaystyle \|d\|_{{\mathcal L}^\infty} \rightarrow 0 \Rightarrow\|d'\|_{{\mathcal L}^\infty} \rightarrow 0. } \label{Dprime}
	\end{eqnarray}
\end{description}

Notice that if system leaves discrete state $q=1$, or $q=2$ then $\eta=h$ and does not depend on $u$ and the condition (\ref{H3c}) is the regularity condition associated with (\ref{sw1a})--(\ref{sw2a}). In \cite{ps2020a} it was shown that the condition (\ref{Dprime}) is equivalent to the condition: there exists $0 < L_d <\infty$ such that
\begin{eqnarray}
    &{\displaystyle  \|d'\|_{{\mathcal L}^\infty} \leq L_d  \|d\|_{{\mathcal L}^\infty}\ \ \forall d\in {\mathcal U}-u.}\label{Dprime_equiv}
\end{eqnarray}
We denote by $N_t(u)$ the number of switching times triggered by a control $u\in {\mathcal U}$ and by $I_t(u)=\{1,2,\ldots,N_t(u)\}$ the set of indices which indicates all switching points. We assume  that for any $u\in {\mathcal U}$ the number of switching points is finite and the last one does not appear at the final time $t_f$

We make the hypothesis 
\begin{description}
		\item[{\bf (H3)}] There exists a nonnegative integer number $I_t^{{\rm max}} < +\infty$ and $\varepsilon_t > 0$ such that
		\begin{eqnarray}
			N_t(u) & \leq & I_t^{{\rm max}},
			\label{Hyp4}\\
			t_t^i & \leq & t_t^{i+1} - \varepsilon_t,\ \forall i\in I_t(u)\ \backslash \{N_t(u)\},\label{Hyp4a} \\
			t_t^{N(u)} & \leq & t_f - \varepsilon_t \label{Hyp4b}
		\end{eqnarray}
		for all $u\in {\mathcal U}$.
\end{description}

\section{Optimal control problem with sliding modes and piecewise smooth controls} 
\label{secOpt}

Taking into account the considerations and definitions presented in the previous section, in particular including the possibility of the system entering (and spending some time) in the sliding mode, we now turn attention to the optimal control problem of interest---{\bf (P)}:
\begin{eqnarray}
&\min_{u\in {\mathcal U}} \phi(x(t_f))& \label{defOptControlProblemCost}
\end{eqnarray}
subject to the constraints
\begin{eqnarray}
\begin{array}{lc}
x' = f_1(x,u) &\ \ {\rm if}\ q=1 \nonumber \\
x' = f_2(x,u) &\ \ {\rm if}\ q=2 \nonumber \\
\left \{
\begin{array}{lll}
x' & = & f_F(x,u) + h_x(x)^Tz \nonumber \\
0 &=& h(x) \nonumber
\end{array}\right .  &\ \ {\rm if}\ q=3 
\end{array}\nonumber \\
\label{SystemEq}
\end{eqnarray}
and the terminal constraints
\begin{eqnarray}
&g^1_i(x(t_f)) = 0\ \forall i\in E & \\
&g^2_j(x(t_f)) \leq 0\ \forall j\in I.& \label{defOptControlProblemEndIneq}
\end{eqnarray}
We also assume that the initial state $x(t_0)=x_0$ is fixed and $T=[t_0,t_f]$.. 

Crucial assumption concerns the set of admissible controls ${\mathcal U}$. Our algorithm requires that control functions are piecewise smooth functions. The justification for that is given in \cite{ps2020a} by referring to properties of Filippov's solutions.

Thus the admissible control $u$ is a function $u$ from the set ${\mathcal U}$, where 
\begin{eqnarray}
&{\displaystyle
	{\mathcal U} =\left \{ u\in {\mathcal C}^N_m[T]:\ u(t)\in
	U\subset  \mathbb{R}^m,\ t\in [t_0,t_f],\right \}}\label{pwc2} \\
&{\displaystyle
	{\mathcal C}_m^N[T] =\left \{ u\in {\mathcal L}^2_m[T]:\ u(t)=\sum_{j=1}^N\psi_j(t)u^N_j(t),\right . }\nonumber \\
&{\displaystyle \left .  
	u_j^N\in \mathcal{C}_m^1[t_{j-1},t_j],\ j=1,\ldots,N,\ t\in[t_0,t_f] \right \}, }
\nonumber
\end{eqnarray}
\begin{eqnarray}
&{\displaystyle 
\psi_j(t)  =   \left \{ \begin{array}{ll}
	1, & {\rm if}\ t\in (t_{j-1},t_j]\nonumber \\
	0, & {\rm if}\ t\not\in (t_{j-1},t_j]\nonumber 
	\end{array} \right .	\ \  \psi_1(t)  =   \left \{ \begin{array}{ll}
	1, & {\rm if}\ t\in [t_0,t_0+(t_f-t_0)/N]\nonumber \\
	0, & {\rm if}\ t\not\in [t_0,t_0+(t_f-t_0)/N]\nonumber 
	\end{array} \right . ,} \nonumber 
\end{eqnarray}
$t_j =  t_0 + j(t_f-t_0)/N,\ j\in \{2,3,\ldots,N\}$.

Regarding the control constraints we make the hypothesis:
\begin{description}		
	\item[{\bf (CC)}] The set of admissible controls ${\mathcal U}$ is a convex compact subset of ${\mathcal L}^2_m[T]$ and there exists $0 < L_d <\infty$ such that
	\begin{eqnarray}
	&{\displaystyle  \|u' - v'[\|_{{\mathcal L}^\infty} \leq L_d  \|u - v\|_{{\mathcal L}^\infty}\ \ \forall u,v\in {\mathcal U}.}\label{Dprime_equivA}
	\end{eqnarray}
\end{description}
For the lack of space we do not analyse in the paper under which conditions imposed on functions $u_j^N\in \mathcal{C}_m^1[t_{j-1},t_j]$---constant, linear, quadratic, etc.---${\bf (CC)}$ is satisfied. It must be stressed that practical application of the presented algorithm requires this form of $u_j^N$ (\cite{ps19a},\cite{ps19b}).

The method we propose for solving the problem ${\bf (P)}$ is based on an exact penalty function. 
By using an exact penalty function approach, instead of solving the problem ${\bf (P)}$ we solve the problem ${\bf (P_{c})}$
\begin{eqnarray}
&{\displaystyle
	\min_{u\in {\mathcal U} } \bar{F}_c(u)}
\label{1e}
\end{eqnarray}
in which the exact penalty function $\bar{F}_c(u)$ is defined as follows
\begin{eqnarray}
&{\displaystyle \bar{F}_c(u) = \bar{F}_0(u) +c \max \left [0,\max_{i\in E} \left |\bar{g}^1_i(u)\right |,\max_{j\in I} \bar{g}^2_j(u)\right ].} \label{2e}
\end{eqnarray}
Here, since under our assumptions $x$ uniquely depends on $u$, we write: $\bar{F}_0(u) = \phi (x^u(t_f))$, $\bar{g}_i^1(u) = g_i^1(x^u(t_f))$, $i\in E$, $\bar{g}_j^2(u) = g_j^2(x^u(t_f))$, $j\in I$.

For fixed $c$ and $u$ the direction finding subproblem, ${\bf P_{c}(u)}$, 
for the problem ${\bf (P_c)}$ is:
\begin{eqnarray}
&{\displaystyle
	\min_{d\in {\mathcal D}_u,\beta\in \mathbb{R}}\left [\left
	\langle \nabla \bar{F}_0(u),d\right \rangle +c \beta + 1/2 
	\| d \|^2_{{\mathcal L}^2} \right ]
} \nonumber
\end{eqnarray}
subject to
\begin{eqnarray}
\left | \bar{g}^1_i(u) +
\left \langle \nabla \bar{g}^1_i(u),d\right \rangle\right |  & \leq &\beta\ \ \forall i\in E
\nonumber \\
\bar{g}^2_j (u) + \left \langle \nabla \bar{g}^2_j(u),d\right \rangle
& \leq & \beta\ \ \forall j\in I.
\nonumber
\end{eqnarray}
Here,
\begin{eqnarray}
&{\displaystyle {\mathcal D}_u \left \{ d\in {\mathcal C}^N_m[T]:\ d\in {\mathcal U}-u\right \}.}\nonumber
\end{eqnarray}

The direction finding subproblem is based on first order approximations to the problem functionals, defined as follows
\begin{eqnarray}
\langle \nabla \bar{F}_0(u),d\rangle & = & \phi_x(x^u(t_f)) y^{x,d}(t_f) \label{Linapp1} \\
\langle \nabla \bar{g}_i^1(u),d\rangle & = & \left (g_i^1\right )_x(x^u(t_f)) y^{x,d}(t_f),\ i\in E \label{Linapp2} \\
\langle \nabla \bar{g}_j^2(u),d\rangle & = & \left (g_j^2\right )_x(x^u(t_f)) y^{x,d}(t_f),\ j\in I. \label{Linapp3}
\end{eqnarray}

The subproblem can be reformulated as an optimization problem 
with the objective function which is strictly convex.
The problem therefore has the unique solution $(\bar{d},\bar{\beta})$.
Since this solution depends on $c$ and $u$, we may define {\it descent function}
$\sigma_{c}(u)$ and {\it penalty test function} $t_{c}(u)$ (see \cite{py99} for details),
to be used to test optimality of a control
$u$ and to adjust $c$, respectively, as 
\begin{eqnarray}
&{\displaystyle
	\sigma_{c}(u) = \left \langle \nabla \bar{F}_0(u),\bar{d}\right \rangle + 
	c\left [\bar{\beta}-M(u)\right ] }\label{sigmaf}
\end{eqnarray}
and                                                          
\begin{eqnarray}
&{\displaystyle
	t_{c}(u) = \sigma_{c}(u) + M(u)/c }\label{testf}
\end{eqnarray}
\label{tesfun2}
for given $c> 0$ and $u\in {\mathcal U}$. Here, 
\begin{eqnarray}
&{\displaystyle
	M(u) = \max\left [0,\max_{i\in E}\left |\bar{g}^1_i(u)\right |,
	\max_{j\in I}\bar{g}^2_j(u)\right ],}\nonumber
\end{eqnarray}

A general algorithm is as follows.
\\[5mm]
{\bf Algorithm}
Fix parameters: $\gamma,\ \eta  \in (0,1)$, $c^0 > 0$, $\kappa >1$.
\begin{enumerate}
	\item Choose the initial control $u_{0}\in {\cal U}$.
	Set $k=0$, $c_{-1}=c^0$.
	\item Let $c_k$ be the smallest number chosen from $\{c_{k-1},\kappa c_{k-1},
	\kappa^2 c_{k-1},\ldots \}$ such that
	the solution $(d_k,\beta_k)$ to the direction finding subproblem  ${\bf P_{c_k}(u_k)}$ 
	satisfies
	\begin{eqnarray}
	&{\displaystyle
		t_{c_k} (u_k) \leq 0. }\label{penadj}
	\end{eqnarray}
	If $\sigma_{c_k}(u_k) =0$ then STOP.
	\item Let $\alpha_k$ be the largest number chosen from the set
	$\{1,\eta,\eta^2,\ldots,\}$ such that $u_{k+1} = u_k + \alpha_k d_k$
	satisfies the relation
	\begin{eqnarray}
	\bar{F}_{c_k}(u_{k+1}) - \bar{F}_{c_k}(u_k)
	& \leq &  \gamma\alpha_k\sigma_{c_k} ( u_k).\nonumber 
	\end{eqnarray}
	Increase k by one. Go to Step 2.
\end{enumerate}

Since
\begin{eqnarray}
&{\displaystyle
	\left \langle \nabla \bar{F}_0(u_k),d\right \rangle + c_k \beta + 1/2 \| d\|^2_{{\mathcal L}^2} \leq
	c_kM(u_k),}\nonumber
\end{eqnarray}
which holds because $0\in {\mathcal U}-u_k$, we also have
\begin{eqnarray}
&{\displaystyle
	\left \langle \nabla \bar{F}_0(u_k),d\right \rangle + c_k\left [\beta -M(u_k)\right ]
	\leq -1/2 \| d\|^2_{{\mathcal L}^2} \leq 0}\label{SigmaPos}
\end{eqnarray}
which implies that the descent function $\sigma_{c_k}(u_k)$ is nonpositive valued at each iteration.

{\it Algorithm} generates a sequence of controls $\{ u_k\}$ and the corresponding sequence of penalty parameters $\{c_k\}$ such that $\{ c_k\}$ is bounded and any accumulation point of
$\{ u_k\}$ satisfies optimality conditions for the problem ${\bf (P_{\bar{c}})}$, i.e. $\sigma_{\bar{c}}(\bar{u})=0$ for the limit point $\bar{u}$ and the limit point of the sequence $\{c_k\}$. 
But $\sigma_{\bar{c}}(\bar{u})=0$ implies that $M(\bar{u})=0$ 
due to the definition (\ref{testf}) and since (\ref{penadj}) holds. 

Since we assume that the set ${\mathcal U}$ is compact we can prove the convergence of {\it Algorithm} without referring to relaxed controls (cf. \cite{py99}). However, we still need to introduce two additional hypotheses. The first one concerns the functions defining the objective and constraints:
\begin{description}
	\item[{\bf (H4)}] $\phi$, $g_i^1$, $i\in E$, $g_j^2$, $j\in I$ are continuously differentiable functions.
\end{description}

The second one is related to a constraint qualification which is needed to prove the convergence of {\it Algorithm}. To this end we first introduce the set
\begin{eqnarray}
&{\displaystyle {\mathcal D} = \left \{ d\in {\mathcal C}^N_m[T]:\ \exists u\in {\mathcal U}\ {\rm such\ that}\ u+d\in {\mathcal U}\right  \} }\nonumber 
\end{eqnarray}
and the set
\begin{eqnarray}
&{\displaystyle  {\mathcal F}(u) = \left \{ d\in {\mathcal D}:\ \max_{j\in I} \left \langle \nabla \bar{g}_j^2(u),d\right \rangle < 0 \right \}. }\nonumber 
\end{eqnarray}
Then the constraint qualification condition takes the form
\begin{description}
\item[{\bf (CQ)}] for each $u\in {\mathcal U}$, ${\mathcal F}(u)\neq \emptyset$, and in the case $E\neq \emptyset$ we have
\begin{eqnarray}
&{\displaystyle 0\in {\rm interior}{\mathcal E}(u)}\label{CQa}
\end{eqnarray}
where
\begin{eqnarray}
&{\displaystyle {\mathcal E}(u) = \left \{ \left \{ \langle \nabla \bar{g}^1_i(u),d\rangle\right \}_{i\in E}\in \mathbb{R}^{|E|}:\ d\in {\mathcal F}(u)\right \}.}\nonumber 
\end{eqnarray}
\end{description}

The constraint qualification ${\bf (CQ)}$ is similar to those stated in \cite{DiPillo} and \cite{Han}, if we observe that we must take into account the constraint $u\in {\mathcal U}$. We will show that under the constraint qualification ${\bf (CQ)}$ any limit point of the sequence $\{u_k\}$ generated by {\it Algorithm}---$\bar{u}\in {\mathcal U}$ (and corresponding to it trajectory $\bar{x}$) will satisfy the conditions
\begin{description}
	\item[{\rm {\bf (NC)}}]:
	\begin{eqnarray}
	&{\displaystyle 0 \leq \min_{d\in {\mathcal D}_{\bar{u}}} \phi_x (\bar{x}(t_f))y^{\bar{x},d}(t_f)}
	\label{NC1a} 
	\end{eqnarray}
	subject to the constraints
	\begin{eqnarray}
	g_i^1(\bar{x}(t_f)) + \left (g_i^1\right )_x(\bar{x}(t_f)) y^{\bar{x},d}(t_f) & = & 0,\ i\in E    
	\label{NC1b} \\
	g_j^2(\bar{x}(t_f)) + \left (g_j^2\right )_x(\bar{x}(t_f)) y^{\bar{x},d}(t_f) & \leq & 0,\ j\in I_{0,\bar{u}}   \label{NC1c}
	\end{eqnarray}
	together with $g_i^1 (\bar{x}(t_f)) = 0$, $i\in E$, $g_j^2 (\bar{x}(t_f))\leq 0$, $j\in I$. Here, 
	\begin{eqnarray}
	&{\displaystyle I_{\varepsilon,u} = \left \{j\in I:\ \bar{g}_j^2(u)\geq \max_{j\in I} \bar{g}_j^2(u) - \varepsilon \right \},\ \varepsilon \geq 0.}\nonumber 
	\end{eqnarray}
\end{description}
In the rest of the paper we say that data for the problem ${\bf (P)}$ satisfy assumptions ${\bf (H1)}$, ${\bf (H2)}$, etc., if at each discrete state functions defining it meet these assumptions. The convergence analysis of {\it Algorithm} is presented in the following theorem.

\newtheorem{t1}{Theorem}[section]
\begin{t1}
	\label{t1}
Assume that data for {\rm {\bf (P)}} satisfy hypotheses ${\bf (H1)}$, ${\bf (H2)}$, ${\bf (H3)}$, ${\bf (H4)}$, ${\bf (CC)}$ and ${\bf (CQ)}$. Let $\{u_k\}$ be a sequence of controls generated by {\it Algorithm} and let $\{c_k\}$ be a sequence of the corresponding penalty parameters. Then 
\begin{enumerate}
	\item [i)] $\{c_k\}$ is a bounded sequence
	\item [ii)] 
		\begin{eqnarray}
		&{\displaystyle \lim_{k\rightarrow\infty}\sigma_{c_k}(u_k)=0,\ \ \lim_{k\rightarrow\infty} M(u_k) = 0.}\label{t1ii}
		\end{eqnarray}
		\item [iii)] Let $\bar{u}$ be the limit point of the sequence $\{u_k\}$ and $\bar{x}$ the trajectory corresponding to $\bar{u}$, then the pair $(\bar{x},\bar{u})$ satisfies ${\bf (NC)}$.
	\end{enumerate}
\end{t1}
The proof of {\it Theorem \ref{t1}} follows the lines of the proof of Theorem 5.1 stated in \cite{pv98}, however for the completeness of presentation it is provided in {\it Appendix \ref{AppendixT}}. 

The Step 2) of {\it Algorithm} is used to approximate the sufficiently large value of the penalty parameter $c$. In particular, if {\bf (CQ)} is satisfied one can prove the following proposition which is needed in order to demonstrate the feasibilities of Steps 2--3 of {\it Algorithm}.

\newtheorem{p1}{Proposition}[section]
\begin{p1}
\label{p1}
	\begin{enumerate}
		\item [(i)]  Assume that data for {\rm {\bf (P)}} satisfy hypotheses ${\bf (H1)}$, ${\bf (H2)}$, ${\bf (H3)}$, ${\bf (H4)}$, ${\bf (CC)}$ and ${\bf (CQ)}$. 
		Then for any $u\in {\mathcal U}$ 
		there exists $\bar{c} >0$ such that for all $c > \bar{c}$
		\begin{eqnarray}
		&{\displaystyle
			t_{c} (u ) \leq 0.}\nonumber
		\end{eqnarray}
		\item [(ii)] Assume that data for {\rm {\bf (P)}} satisfy hypotheses ${\bf (H1)}$, ${\bf (H2)}$, ${\bf (H3)}$ and ${\bf (H4)}$. Then
		for any $u\in {\mathcal U}$
		and $c >0$ such that $\sigma_c (u) < 0$ there exists $\bar{\alpha} > 0$
		such that if $\alpha\in [0,\bar{\alpha})$ then
		\begin{eqnarray}
		\tilde{F}_c (\tilde{u}) - \tilde{F}_c (u) & \leq & \gamma \alpha \sigma_c (u)\nonumber
		\end{eqnarray}
		where 
		$\tilde{u} = u +\alpha d$, $(d, \beta)$ is the solution to the direction finding subproblem corresponding to $c$ and $u$ and $\sigma_c(u)$ is defined by (\ref{sigmaf}).
	\end{enumerate}
\end{p1}
Since for $\bar{\alpha}\leq 1$ $\tilde{u}\in {\mathcal U}$ part {\it (ii)} of the proposition states that $\sigma_c(u)=0$ ($(d\equiv 0,\beta=0))$ are the necessary optimality conditions for the problem ${\bf (P_c)}$. 

The proof of {\it Proposition \ref{p1}} is given in {\it Appendix \ref{AppendixP}}. 

\section{Optimality results for the exact penalty function}
\label{exactPenalty}

{\it Algorithm} finds a control $\bar{u}$ which satisfies the necessary optimality conditions ${\bf (NC)}$. In order to show that a control fulfilling these conditions is a good candidate for a local solution of the problem ${\bf (P)}$ we will present analysis which shows relations between local solutions of problems ${\bf (P_c)}$ and the problem ${\bf (P)}$. 
We hope that the presented analysis emphasizes further the importance of the condition ${\bf (CQ)}$. 

The proofs of {\it Theorem \ref{t1}} and {\it Proposition \ref{p1}} heavily use the implications of the constraint qualification {\bf (CQ)} stated in the following lemma. The role of {\bf (CQ)} is to ensure uniform boundedness of penalty parameter values.
The proof of the lemma closely follows the lines of the proof of Lemma 6.1 stated in \cite{pv98}. However due to the importance of the lemma we present its proof. Furthermore, we will use parts of the proof to establish the results stated in {\it Theorem \ref{t1p}}.

\newtheorem{l2}{Lemma}[section]
\begin{l2}
\label{l2}
	Assume ${\bf (H1)}$, ${\bf (H2)}$, ${\bf (H3)}$, ${\bf (H4)}$, ${\bf (CC)}$ and ${\bf (CQ)}$. 
	For any control $\tilde{u}\in{\mathcal U}$
	there exist a neighborhood ${\mathcal O}(\tilde{u})$ of $\tilde{u}$, 
	$K_1 > 0$ and $K_2 > 0$ with
	the following properties: given any 
	$u\in {\mathcal U}$ 
	such
	that
	$u\in{\mathcal O}(\tilde{u})$
	there exists $v\in {\mathcal U}$ such that
	\begin{eqnarray}
	&{\displaystyle \max_{i\in E} \left | \bar{g}^1_i(u) + \left \langle \nabla \bar{g}^1_i(u),
		v- u\right \rangle \right | - \max_{i\in
			E}
		\left |\bar{g}^1_i (u) \right | \leq -  K_1 M(u) }
	\label{l532}\\
	&{\displaystyle 
		\max_{j\in I} 
		\left \langle \nabla \bar{g}^2_j(u),v-u \right \rangle \leq  - K_1 M(u) }
	\label{l533} \\
	&{\displaystyle 
		{\rm and}\ \ \left \| v -u \right \|_{{\mathcal L}^\infty} \leq \ \ K_2 M(u).}
	\label{l535}
	\end{eqnarray}
\end{l2}
\begin{proof} 
Take any $u\in {\cal U}$.
	Let $r>0$ be a number such that
	\begin{eqnarray}
	&{\displaystyle
		M(u) < r}\nonumber
	\end{eqnarray}
	for all $u\in {\cal U}$.
	We deduce from ${\bf (CQ)}$ that there is a simplex
	in ${\cal E}(u)\subset
	{\cal R}^{n_E}$ ($n_E= |E|$) with vertices $\{e_j\}_{j=0}^{n_E}$ which contains $0$
	as an interior point. By definition of ${\cal E}(u)$, there exist
	$d_0,\ldots,d_{n_E}\in {\cal D}$ and $\delta > 0$ such that for $j=0,\ldots,n_E$
	\begin{eqnarray}
	&{\displaystyle 
	\left \{\left \langle \nabla \bar{g}^1_i(u),
	d_j\right \rangle\right \}_{i\in E} =  e_j,\ \ 
	\max_{i\in I} \left \langle \nabla \bar{g}^2_i(u),
	d_j \right \rangle  \leq  - \delta.}
	\nonumber 
	\end{eqnarray}
	
	Let $(\lambda_0,\lambda_1,\ldots,\lambda_{n_E})$ be the barycentric coordinates of $0$
	w.r.t. the vertices $e_j$ of the simplex, i.e.
	\begin{eqnarray}
	&{\displaystyle
		0\ \left (=\sum_{j=0}^{n_E} \lambda_j e_j\right ) = \nabla
		\bar{g}^1 (u)\circ \sum_{j=0}^{n_E} \lambda_j d_j .}\nonumber
	\end{eqnarray}
	Here
	\begin{eqnarray}
	\nabla \bar{g}^1(u)\circ d & := & \left \{\left \langle \nabla \bar{g}^1_i(u),d\right \rangle
	\right \}_{i\in
		E}.\nonumber
	\end{eqnarray}
	We shall also write
	\begin{eqnarray}
	\bar{g}^1 (u) & := & \left \{ \bar{g}^1_i(u)\right \}_{i\in E}.\nonumber
	\end{eqnarray}
	
	Since the vertices are in general position and $0$ is an interior point, the
	$\lambda_i$'s are all positive and we may find $\delta_1 > 0$ such that
	\begin{eqnarray}
	&{\displaystyle
		\left (\lambda_0 - \sum_{j=1}^{n_E} \alpha_j,\lambda_1 + \alpha_1,\ldots,
		\lambda_{n_E} + \alpha_{n_E} \right )}
	\nonumber \\
	&{\displaystyle
		\in \left \{ \gamma\in {\cal R}^{n_E+1}:\ \gamma_j \geq 0\ \forall
		j,\ \sum_{j=0}^{n_E}\gamma_j =1 \right \}} \nonumber
	\end{eqnarray}
	whenever $\alpha\in {\cal B}(0,\delta_1)\subset {\cal R}^{n_E}$. (${\cal B}(0,\delta_1)$
	is a ball with radius $\delta_1$.)
	Furthermore, the $n_E\times n_E$ matrix $P(u)$
	defined by
	\begin{eqnarray}
	&{\displaystyle
		P(u)\alpha := \sum_{j=1}^{n_E} \nabla \bar{g}^1(u)\circ \alpha_j (d_j-d_0) }
	\label{b1}
	\end{eqnarray}
	is invertible for $u=\tilde{u}$, from the definition of $d_j,\ j= 1,\ldots,n_E$.
	
	In consequence of hypothesis ${\bf (CQ)}$ and in view of the continuity properties of
	the mapping $u\rightarrow y^{x,d}$ for fixed $d$ 
	we may choose a neighborhood ${\cal O}(\tilde{u})$ of $\tilde{u}$ in ${\cal U}$ and
	numbers 
	$\hat{r}\geq r$ and 
    $\delta_2\in (0,\hat{r}]$ such that for any
	$u\in {\cal U}$ satisfying $u\in {\cal O}(\tilde{u})$
	\begin{eqnarray}
	(i) \hspace{5mm}& & \max_{i\in I}\left \langle \nabla \bar{g}^2_i(u),v_j - u\right \rangle \leq
	- \delta/2\ \ \forall j   \nonumber\\
	(ii)\hspace{5mm} & & P(u)\ \ {\rm is\ invertible} \nonumber\\
	(iii)\hspace{5mm} & & \left \| P(u)^{-1} \nabla \bar{g}^1(u)\circ
	\left (\left (\sum_{j=0}^{n_E} \lambda_j v_j\right )-u\right )
	\right \| \leq \delta_1/2\nonumber\\
	(iv)\hspace{5mm} & & \delta_2 \left \| P(u)^{-1} 
	\right \| n_E^{1/2} \leq \delta_1/2 .\nonumber
	\end{eqnarray}
	(In {\it (iv)} the norm is the Frobenius norm.
	Here the controls $v_j\in {\cal U}$, $j=0,\ldots,n_E$ are defined to be
	\begin{eqnarray}
	&{\displaystyle
        v_j:= u + d_j(u).}\nonumber
	\end{eqnarray}
	
	Now suppose that the control $u$ is not feasible.
	Set
	\begin{eqnarray}
	&{\displaystyle
		\alpha = P(u)^{-1} \left [ - \nabla \bar{g}^1(u)\circ
		\left (\left (\sum_{j=0}^{n_E}\lambda_j v_j\right ) -u\right )
		- \delta_2 M(u)^{-1}
		\bar{g}^1(u) \right ]. }\label{b2}
	\end{eqnarray}
	Notice that, by properties {\it (iii)} and {\it (iv)},
	$\| \alpha \| \leq \delta_1$. Also set
	\begin{eqnarray}
	&{\displaystyle
		\hat{v} = v_0 + \sum_{j=1}^{n_E} (\lambda_j + \alpha_j)(v_j-v_0).}\label{T1p_p1}
	\end{eqnarray}
	Because $\| \alpha \| \leq \delta_1$ we have that $\hat{v}\in {\cal U}$.
	Finally we define $v$ to be
	\begin{eqnarray}
	&{\displaystyle
		v = u + \left (M(u)/\hat{r}\right )(\hat{v}-u).}\nonumber
	\end{eqnarray}
	Since $M(u)/\hat{r} \leq 1$, it follows that
	$v\in {\cal U}$. We now verify that this control function has the required
	properties.                     
	
	Notice first that
	\begin{eqnarray}
	&{\displaystyle
		\left \| v - u \right \|_{{\mathcal L}^\infty} \leq (2d/\hat{r})M(u),}\label{b3}
	\end{eqnarray}
	where $d$ is a bound on the norms of elements in ${\cal U}$.
	
	We have from (\ref{b1}) and (\ref{b2}) that
	\begin{eqnarray}
	P(u)\alpha &  = & \nabla\bar{g}^1(u)\circ \sum_{j=1}^{n_E} 
	\alpha_j(v_j- v_0) 
	\nonumber \\
	& = & - \nabla \bar{g}^1(u)\circ \left (\sum_{j=1}^{n_E} \lambda_j \left (v_j-v_0\right ) 
	+ v_0 -
	u\right ) - \delta_2 M(u)^{-1}
	\bar{g}^1(u).\nonumber
	\end{eqnarray}
	By definition of $\hat{v}$,
	\begin{eqnarray}
	&{\displaystyle
		\nabla \bar{g}^1(u)\circ (\hat{v}-u) = -\delta_2 M(u)^{-1}\bar{g}^1(u).}\label{T1p_p2}
	\end{eqnarray}
	But then
	\begin{eqnarray}
	&{\displaystyle
		\nabla \bar{g}^1(u)\circ (v - u) = -(\delta_2/
		\hat{r})\bar{g}^1(u).}
	\nonumber
	\end{eqnarray}
	Since $\delta_2/\hat{r}\leq 1$, it follows that
	\begin{eqnarray}
	&{\displaystyle
		\max_{i\in E} \left | \bar{g}^1_i(u) + \left \langle \nabla \bar{g}^1_i(u),
		v-u\right \rangle \right |
		- \max_{i\in E} \left | \bar{g}^1_i(u)\right |
		\leq - (\delta_2/\hat{r})M(u).}\label{b4}
	\end{eqnarray}
	
	We deduce from property {\it (i)} that
	\begin{eqnarray}
	&{\displaystyle
	\left \langle \nabla \bar{g}^2_j(u),v_0 + \sum_{i=1}^{n_E}\left (\lambda_i +
	\alpha_i\right )\left (v_i-v_0\right )-u\right \rangle 
    \leq - \delta/2,\ \forall j\in I.}
	\label{T1p_p3}
	\end{eqnarray}
Since $M(u)/\hat{r}\leq 1$ 
	we deduce that
	\begin{eqnarray}
	&{\displaystyle
		\min\left [0,\bar{g}^2_j(u)\right ]+
		\left \langle \nabla \bar{g}^2_j(u),v-u\right \rangle \leq 
		-\left (\delta/\left (2\hat{r}\right )\right )M(u)\ \forall j\in I.}\label{b5}
	\end{eqnarray}
	
	Surveying inequalities (\ref{b3})--(\ref{b5}),
	we see that $v$ satisfies all relevant
	conditions for completion of the proof, 
	when we set $K_1 = \min \{\delta_2/\hat{r},\delta/(2\hat{r})\}$
	and $K_2 = 2d/\hat{r}$,
	numbers whose magnitudes do not depend on our choice of $u$.
\end{proof}


First we show, using {\it Lemma \ref{l2}}, that if $\bar{u}$ is a strict local solution to the problem ${\bf (P)}$ then for sufficiently large values of $c$ it is a local solution to the problem ${\bf (P_c)}$.

\newtheorem{t1p}{Theorem}[section]
\begin{t1p}
\label{t1p}
Assume that data for {\rm {\bf (P)}} satisfy hypotheses ${\bf (H1)}$, ${\bf (H2)}$, ${\bf (H3)}$, ${\bf (H4)}$ and ${\bf (CC)}$. Suppose that $\bar{u}$ is a strict local minimum of the problem ${\bf (P)}$ (which means that $\bar{u}$ is feasible with respect to all constraints, in particular $\bar{u}\in {\mathcal U}$) and that at the point $\bar{u}$ constraint qualification ${\bf (CQ)}$ holds. Then there exists an $\bar{c} > 0$ such that for all $c\geq \bar{c}$, $\bar{u}$ is a local minimum of ${\bf (P_c)}$ on  ${\cal U}$. 
\end{t1p}
\begin{proof}
	
For the simplicity of presentation we assume that $I = \emptyset$. The proof goes along the lines of the proof of Theorem 4.4 in \cite{Han}.
	
Suppose that $\bar{u}$ is the strict local minimum of the problem ${\bf (P)}$ so there exists the open ball with the centre $\bar{u}$ and  the radius $\bar{\varepsilon} > 0$---${\mathcal B}(\bar{u},\bar{\varepsilon})$ such that $\bar{F}_0(u) > \bar{F}_0 (\bar{u})$ for all $u\in {\mathcal B}(\bar{u},\bar{\varepsilon}) \cap {\mathcal U}\cap {\mathcal G}_E$, $u\neq \bar{u}$, where ${\mathcal G}_E = \{ u\in {\mathcal C}^N_m[T]:\ \bar{g}^1_i(u) =0,\ i\in E\}$. According to Theorem in \cite{Piet} for sufficiently large values of $c$ there exist $\varepsilon (c) > 0$ and $u(c)$ such that $u(c)$ is a local minimum of the problem ${\bf (P_c)}$ on the set ${\mathcal B}(\bar{u},\varepsilon (c))\cap {\mathcal U}$. Moreover it follows from the theorem that $\varepsilon (c) \rightarrow_{c\rightarrow \infty} 0$ so for $c$ sufficiently large we will have $\varepsilon (c) \leq \bar{\varepsilon}$. Suppose that $u(c)$ is feasible with respect to the constraints $\bar{g}^1_i(u) = 0,\ i\in E$. Then we will have
	\begin{eqnarray}
	&{\displaystyle \bar{F}_c(u(c)) = \bar{F}_0(u(c)) \leq \bar{F}_c(\bar{u}),} \nonumber
	\end{eqnarray}
	which means that we must have $u(c) = \bar{u}$ since $\bar{u}$ is a strict local minimizer on the set ${\mathcal B}(\bar{u},\varepsilon (c))\cap {\mathcal U}\cap {\mathcal G}_E$ so it is the only minimizer in the set. This means $\bar{u}$ is the local minimizer for the problem ${\bf (P_c)}$ on the set ${\mathcal B}(\bar{u},\varepsilon (c))\cap {\mathcal U}$.
	
	Therefore, we have to show that for sufficiently large values of $c$ vectors $u(c)$ are feasible with respect to the constraints $\bar{g}^1_i(u) = 0,\ i\in E$. Suppose that it is not true, thus for any $c\rightarrow \infty$ there exists an index $i_c\in E$ such that $\bar{g}^1_{i_c}(u(c)) \neq 0$.  Since at $\bar{u}$ ${\bf (CQ)}$ holds there exists a neighborhood ${\mathcal B}(\bar{u},\varepsilon_1)$, $\varepsilon_1 > 0$ such that for any $u\in {\mathcal B}(\bar{u},\varepsilon_1)$ there exists $v\in {\mathcal U}$ and $L_1 > 0$ with the property (We take as $v$ the control $\hat{v}\in {\mathcal U}$ introduced in the proof of {\it Lemmma \ref{l2}}. For the control $\hat{v}$ defined by (\ref{T1p_p1}) the equalities (\ref{T1p_p2}) and the inequalities (\ref{T1p_p3}) hold so we can take $L_1 = \delta_2$.)
	\begin{eqnarray}
	&{\displaystyle \left \langle \nabla \bar{g}^1_i(u),v- u\right \rangle   =  -L_1 \bar{g}_i^1(u)/M(u),\ i\in E.      }\nonumber 
	\end{eqnarray}
	
	Suppose that $\bar{i}_c(u) \in E$ is such that $\left |\bar{g}^1_{\bar{i}_c(u)} (u) \right | = M(u)$ (notice that we have assumed that $I = \emptyset$). Then we have
	\begin{eqnarray}
	&{\displaystyle  \left \langle \nabla \bar{g}^1_i(u),v- u\right \rangle = \left \{ 
		\begin{array}{ll}
		-L_1 & {\rm if}\ i=\bar{i}_c(u)\ \ {\rm and}\ \bar{g}^1_{\bar{i}_c}(u) = M(u) \nonumber \\
		L_1   & {\rm if}\ i=\bar{i}_c(u)\ \ {\rm and}\ \bar{g}^1_{\bar{i}_c}(u) = - M(u) \nonumber \\                   
		-L_1 \bar{g}^1_i(u)/M(u) & {\rm if}\ \bar{g}^1_i(u) \geq 0 \nonumber \\
		L_1 \bar{g}^1_i(u)/M(u) & {\rm if}\ \bar{g}^1_i(u) < 0 \nonumber
		\end{array} 
		\right. .} \label{ep1}                
	\end{eqnarray}
	
For sufficiently large $c$ $u(c)\in {\mathcal B}(\bar{u},\varepsilon_1)$ and $u(c)$ is not feasible with respect to constraints $\bar{g}^1_i(u) = 0,\ i\in E$. From {\it Lemma \ref{l2}} for any $c$ there exists $v \in {\mathcal U}$ for which the directional derivative of $\bar{F}_c(u)$ in the direction $v-u$ can be evaluated according to the formula (Danskin's formula, see, for example, Theorem 10.22 in \cite{cla2013})
	
\begin{eqnarray}
	&{\displaystyle 
		D\bar{F}_c(u;v-u) = \left \langle \nabla \bar{F}_0(u),v-u \right \rangle + c \max \left [  \max_{i \in I_{+}(u)} \left \langle \nabla \bar{g}^1_i(u),v-u \right \rangle ,\right . }\nonumber \\
	&{\displaystyle \left . \max_{i \in I_{-}(u)} \left [ - \left \langle \nabla \bar{g}^1_i(u),v-u \right \rangle \right ] \right ], }\nonumber 
	\end{eqnarray}
	where 
	\begin{eqnarray}
 I_{+}(u) &  = & \left \{ i \in I(u):\ \bar{g}_i^1(u)  = M(u)\right \}\nonumber \\
 I_{-}(u) &  = & \left \{ i \in I(u):\  \bar{g}_i^1(u) = - M(u)\right \}\nonumber 
	\end{eqnarray}
	(notice that $I_{+}(u) \cap I_{-}(u) = \emptyset$ since $M(u) > 0$).
	
	From (\ref{ep1}) it follows that
	\begin{eqnarray}
	&{\displaystyle D\bar{F}_c(u;v-u) = \left \langle \nabla \bar{F}_0(u),v-u \right \rangle - K_1 c \leq  \|\bar{F}_0(u) \| \|v-u \| - L_1c. }\label{ep2}
	\end{eqnarray}
 
	Therefore, for sufficiently large values of $c$, since $v\in {\mathcal U}$ which is bounded,  we have that $u(c)\in {\mathcal B}(\bar{u},\varepsilon_1)$ and at the same time, from (\ref{ep2}),
	\begin{eqnarray}
	&{\displaystyle D\bar{F}_c(u(c);v-u(c)) < 0} \label{ep3}
	\end{eqnarray}
	and, since $v\in {\mathcal U}$ (which is convex), (\ref{ep3}) contradicts the assumption that $u(c)$ is a local minimizer of the problem ${\bf (P_c)}$. It means that for sufficiently large values of $c$ $u(c)$ is feasible with respect to the constraints $\bar{g}^1_i(u) = 0,\ i\in E$.
\end{proof}

Using parts of the proof of {\it Theorem \ref{t1p}} we are able to prove the following theorem (cf. Theorem 4.1 in \cite{DiPillo}).
\newtheorem{t2p}[t1p]{Theorem}
\begin{t2p}
\label{t2p}
Assume that data for {\rm {\bf (P)}} satisfy hypotheses ${\bf (H1)}$, ${\bf (H2)}$, ${\bf (H3)}$, ${\bf (H4)}$, ${\bf (CC)}$ and ${\bf (CQ)}$. Then there exists $\bar{c} > 0$ such that for $c\geq \bar{c}$ if $u(c)$ is a local minimum point of the problem ${\bf (P_c)}$ then $u(c)$ is also a local minimum point of the problem ${\bf (P)}$. 
\end{t2p}
\begin{proof}
Using arguments similar to those applied in the proof of {\it Theorem \ref{t1p}}, and also those which are presented in the proof of Proposition 3.3 stated in \cite{DiPillo}, we can show that there exists $\bar{c} > 0$ such that for any $c \geq \bar{c}$, if $u(c)$ is a local solution to the problem ${\bf (P_c)}$ then $u(c)$ is feasible with respect to the constraints of the problem ${\bf (P)}$. It means that there exists ${\mathcal B}(u(c),\varepsilon)$ ($\varepsilon > 0$) such that
\begin{eqnarray}
&{\displaystyle \bar{F}_0(u(c)) = \bar{F}_c(u(c)) \leq \bar{F}_c(u),\ \forall u\in {\mathcal B}(u(c),\varepsilon)\cap {\cal U},\ c\geq \bar{c},}\nonumber
\end{eqnarray}
which implies that
\begin{eqnarray}
&{\displaystyle \bar{F}_0(u(c)) = \bar{F}_c(u(c)) \leq  \bar{F}_c(u) = \bar{F}_0(u),\ \forall u\in {\mathcal B}(u(c),\varepsilon)\cap {\mathcal U}\cap {\mathcal G},} \label{ep4}
\end{eqnarray}
$c\geq \bar{c}$, where ${\mathcal G} = \{u\in {\mathcal C}^N_m[T]:\ \bar{g}^1_i(u) = 0,\ i\in E,\ \bar{g}^2_j(u) \leq 0,\ j\in I\}$. But (\ref{ep4}) states that $u(c)$ is a strict local minimum for the problem ${\bf (P)}$.
\end{proof} 

The aim of {\it Algorithm} is to indicate a candidate for a local solution of the problem ${\bf (P)}$. It finds a control ${\bar{u}}$ which satisfies the optimality conditions ${\bf (NC)}$, in particular it is feasible with respect to the constraints of the problem ${\bf (P)}$. {\it Algorithm} generates the sequence  $\{c_k\}$ which has the limit point (since $\{c_k\}$ is bounded and non--decreasing) which we denote by $\bar{c}$. 

 
If we assume that $\bar{u}$ is a local solution for the problem $\min_{u\in {\mathcal U}}\bar{F}_c(u)$, then from {\it Theorem \ref{t2p}} it follows that $\bar{u}$ is a local solution to the problem ${\bf (P)}$. Therefore, {\it Algorithm} designates a candidate for a local solution of the problem ${\bf (P)}$.

%% file: sec_adjointEqsSC_P2.tex
\section{Adjoint equations}
\label{secAdjoint}
The efficient implementation of {\it Algorithm} needs computationally tractable algorithm for solving the direction finding subproblem ${\bf P_c(u)}$ which at every iteration requires the evaluation of $\langle \nabla \bar{F}_0(u),d\rangle$, $\langle \nabla \bar{g}_i^1(u),d\rangle$, $i\in E$,  $\langle \nabla \bar{g}_j^2(u),d\rangle$, $j\in I$. If we do that by using solutions to the linearized equations then each update of $d$ would request solving these equations. If we use adjoint equations the evaluation comes down to the calculations of scalar products since $\nabla \bar{F}_0(u)$, $\nabla \bar{g}_i^1(u)$, $i\in E$, $\nabla \bar{g}_j^2(u)$, $j\in I$ are available prior to solving ${\bf P_c(u)}$. It is sufficient to show how $\nabla \bar{F}_0(u)$ could be calculated with the help of adjoint equations since the adjoint equations for the other functionals will differ only by terminal conditions.

We formulate the adjoint equations for three cases, they are derived on the basis of a variational approach presented in \cite{bh1975} and by taking into account our results on solutions to linearized equations of ordinary differential equations and differential--algebraic equations of the special form arising when a system remains in the sliding mode (\cite{ps2020a}).  

{\it Case 1--3).} In the first case we assume that in the time interval $ [t_0,t_t] $ the system evolves according to 
\begin{equation}
x' = f_1(x,u).
\label{eqAdjEqsFirstCaseOde}
\end{equation}
At a transition time $ t_t $ the continuous state trajectory meets the switching surface such that the following condition holds
\begin{equation}
h(x(t_t)) = 0. 
\label{eqAdjEqsFirstCaseTrans}
\end{equation} 
After the transition the system evolves according to DAEs 
\begin{eqnarray}
x' &=& f_F(x,u) + h_x^T(x)z \label{eqAdjEqsFirstCaseDaeDiff} \\
0 &=& h(x) \label{eqAdjEqsFirstCaseDaeAlg}
\end{eqnarray}
up to an ending time $ t_f $. 

To derive the adjoint equations we construct the following augmented functional
\begin{eqnarray}
&{\displaystyle \Phi(x,z,u,\lambda_f,\lambda_h,\pi) = \phi(x(t_f)) + \pi h(x(t_t))  + } \nonumber \\
&{\displaystyle \int_{t_0}^{t_t^-} \lambda_f^T(t) \left( x'(t) - f_1(x(t),u(t)) \right) dt +}\nonumber   \\
&{\displaystyle \int_{t_t^+}^{t_f}\left [ \lambda_f^T(t) \left( x'(t) - f_F(x(t),u(t)) - h_x^T(x(t))z(t) \right)  
+\lambda_h^T(t) h(x(t))\right ] dt.} \nonumber
\end{eqnarray} 
Suppose now that we calculate the variation of the augmented functional with respect to $(x,z)$ only.
\begin{eqnarray}
&{\displaystyle d\Phi(x,z,u,\lambda_f,\lambda_h,\pi) = \phi_x(x(t_f))y^{x,d}(t_f) + \pi h_x(x(t_t))dx(t_t) + }\nonumber \\ 
&{\displaystyle \lambda_f^T(t_t^-) \left( x'(t_t^-)
- f_1(x(t_t^-),u(t_t^-)) \right) dt_t }\nonumber \\
&{\displaystyle 
 +d\left [\int_{t_0}^{t_t^-} \lambda_f^T(t) \left( x'(t) - f_1(x(t),u(t)) \right) dt\right ] }\nonumber \\
&{\displaystyle -\lambda_f^T(t_t^+) \left( x'(t_t^+) - f_F(x(t_t^+),u(t_t^+)) - \right .}\nonumber \\
&{\displaystyle \left . h_x^T(x(t_t^+))z(t_t^+) \right)dt_t   -\lambda_h^T(t_t^+) h(x(t_t^+)) dt_t }\nonumber \\
&{\displaystyle + d \left [\int_{t_t^+}^{t_f}\left [ \lambda_f^T(t) \left( x'(t) - f_F(x(t),u(t)) -\right . \right . \right .}\nonumber \\
&{\displaystyle \left . \left . \left .h_x^T(x(t))z(t) \right) 
+\lambda_h^T(t) h(x(t))\right ] dt\right ].} \nonumber
\end{eqnarray} 
By integrating by parts the formulas $ \int \lambda_f(t)x(t)dt $ (since $\lambda_f$ is sufficiently smooth) we obtain 
\begin{eqnarray}
&{\displaystyle d\Phi(x,z,u,\lambda_f,\lambda_h,\pi) = \phi_x(x(t_f))y^{x,d}(t_f) + \pi h_x(x(t_t))dx(t_t)}\nonumber \\
&{\displaystyle  +\lambda_f^T(t_t^-) \left( x'(t_t^-) - f_1(x(t_t^-),u(t_t^-)) \right) dt_t + d\left [ \left[\lambda_f^T(t)x(t)\right]_{t_0}^{t_t^-}\right ] 
-}\nonumber \\
&{\displaystyle d\left [ \int_{t_0}^{t_t^-}\left [(\lambda_f^T)'(t)  x(t) + \lambda_f^T(t) f_1(x(t),u(t))\right ] dt \right ] }\nonumber \\
&{\displaystyle -\lambda_f^T(t_t^+) \left( x'(t_t^+) - f_F(x(t_t^+),u(t_t^+))  - h_x^T(x(t_t^+))z(t_t^+) \right)dt_t }\nonumber\\
&{\displaystyle -\lambda_h^T(t_t^+) h(x(t_t^+)) dt_t + d \left [ \left[\lambda_f^T(t)x(t)\right]_{t_t^+}^{t_f} \right ] - }\nonumber \\
&{\displaystyle d\left [ \int_{t_t^+}^{t_f} \left [ (\lambda_f^T)'(t)  x(t)  + \lambda_f^T(t) \left( f_F(x(t),u(t)) +
\right .\right .\right .}\nonumber \\
&{\displaystyle \left . \left . \left . 
h_x^T(x(t))z(t) \right)  - \lambda_h^T(t) h(x(t))\right ] dt\right ].} \nonumber
\end{eqnarray} 
Expanding further the variations (also with respect to u) and taking into account the initial conditions of the linearized equations we obtain
\begin{eqnarray}
&{\displaystyle d\Phi(x,z,u,\lambda_f,\lambda_h,\pi) = \phi_x(x(t_f))y^{x,d}(t_f) + \pi h_x(x(t_t))dx(t_t) +}\nonumber \\
&{\displaystyle \lambda_f^T(t_t^-)x'(t_t^-)dt_t - \lambda_f^T(t_t^-)f_1(x(t_t^-),u(t_t^-))dt_t  +\lambda_f^T(t_t^-)y^{x,d}(t_t^-) }\nonumber \\
&{\displaystyle -\int_{t_0}^{t_t^-}\left [ (\lambda_f^T)'(t) y^{x,d}(t) + \lambda_f^T(t) (f_1)_x(x(t),u(t))y^{x,d}(t) \right .}\nonumber \\
&{\displaystyle \left . +\lambda_f^T(t) (f_1)_u(x(t),u(t))d(t) \right ]dt -\lambda_f^T(t_t^+)x'(t_t^+)dt_t +}\nonumber \\
&{\displaystyle \lambda_f^T(t_t^+)f_F(x(t_t^+),u(t_t^+))dt_t +\lambda_f^T(t_t^+)h_x^T(x(t_t^+))z(t_t^+)dt_t -}\nonumber \\
&{\displaystyle  \lambda_h^T(t_t^+)h(x(t_t^+))dt_t + \lambda_f^T(t_f)y^{x,d}(t_f) - \lambda_f^T(t_t^+)y^{x,d}(t_t^+) - }\nonumber\\
&{\displaystyle \int_{t_t^+}^{t_f} \left [(\lambda_f^T)'(t) y^{x,d}(t) + \lambda_f^T(t) (f_F)_x(x(t),u(t)) y^{x,d}(t) \right . }\nonumber \\
&{\displaystyle + \lambda_f^T(t) (f_F)_u(x(t),u(t))d(t) + \lambda_f^T(t) h_x^T(x(t))y^{z,d}(t) +}\nonumber \\
&{\displaystyle \left . \lambda_f^T(t) (h_x^T(x(t))z(t))_xy^{x,d}(t)-\lambda_h^T(t) h_x(x(t))y^{x,d}(t) \right ]dt.} \nonumber
\end{eqnarray}  

Now, if we take into account the formula (\ref{Lin5})--(\ref{Lin6}) for the differential of $ x(t_t) $,
and rearrange the components with respect to differentials $ dx(t_t),\ dt_t $ and variations $ y^{x,d}(t_f),\ y^{x,d}(t),\ y^{z,d}(t),\ d(t) $ we come to
\begin{eqnarray}
&{\displaystyle d\Phi(x,z,u,\lambda_f,\lambda_h,\pi) = \left( \phi_x(x(t_f)) +\lambda_f^T(t_f) \right)y^{x,d}(t_f) + \left( \pi h_x(x(t_t)) + \right .}\nonumber \\
&{\displaystyle \left . \lambda_f^T(t_t^-) -\lambda_f^T(t_t^+) \right) dx(t_t) + \left( -\lambda_f^T(t_t^-)f_1(x(t_t^-),u(t_t^-)) + \right .} \nonumber \\ 
&{\displaystyle \left . \lambda_f^T(t_t^+)f_F(x(t_t^+),u(t_t^+))  + \lambda_f^T(t_t^+)h_x^T(x(t_t^+))z(t_t^+) - \lambda_h^T(t_t^+)h(x(t_t^+)) \right) dt_t }\nonumber \\
&{\displaystyle -\int_{t_0}^{t_t^-} \left ( (\lambda_f^T)'(t) + \lambda_f^T(t) (f_1)_x(x(t),u(t)) y^{x,d}(t)  + \right .}\nonumber \\
&{\displaystyle \left . \lambda_f^T(t) (f_1)_u(x(t),u(t))d(t) \right )dt }\nonumber \\
&{\displaystyle -\int_{t_t^+}^{t_f} \left [ \left( (\lambda_f^T)'(t) +  \lambda_f^T(t) (f_F)_x(x(t),u(t)) + \lambda_f^T(t)  (h_x^T(x(t))z(t))_x -\right . \right .} \nonumber \\
&{\displaystyle \left . \left . \lambda_h^T(t) h_x(x(t)) \right) y^{x,d}(t) + \lambda_f^T(t) (f_F)_u(x(t),u(t))d(t) + \right . }\nonumber \\
&{\displaystyle \left . \lambda_f^T(t) h_x^T(x(t))y^{z,d}(t)\right ]  dt.} \nonumber  
\end{eqnarray} 

We can now state conditions for adjoint equations such that the expressions with differentials $ dx(t_t),\ dt_t $ and variations $ y^{x,d}(t_f),\ y^{x,d}(t),\ y^{z,d}(t) $ disappear and only the coefficients with variation $ d(t) $ remain. These components will state a formula to efficiently calculate the gradient of $ \phi(x(t_f)) $ with respect to controls.

Let us start with components that are present under the integrals. We want to zero the following components
\begin{eqnarray}
&{\displaystyle \left( (\lambda_f^T)'(t) + \lambda_f^T(t) (f_1)_x(x(t),u(t)) \right) y^{x,d}(t),\ t\in[t_0,t_t^-], }\nonumber \\
&{\displaystyle \left( (\lambda_f^T)'(t) + \lambda_f^T(t) (f_F)_x(x(t),u(t)) + \right .}\nonumber \\
&{\displaystyle \left . \lambda_f^T(t) (h_x^T(x(t))z(t))_x  -\lambda_h^T(t) h_x(x(t)) \right) y^{x,d}(t),\ t\in[t_t^+,t_f],}\nonumber \\
&{\displaystyle \lambda_f^T(t) h_x^T(x(t))y^{z,d}(t),\ t\in[t_t^+,t_f] .}\nonumber 
\end{eqnarray}
We achieve that by assuming that
\begin{equation}
(\lambda_f^T)'(t) = - \lambda_f^T(t) (f_1)_x(x(t),u(t)),\ t\in[t_0,t_t^-]  
\label{eqAdjEqsFirstCaseAdjOde}
\end{equation}
and 
\begin{eqnarray}
(\lambda_f^T)'(t) &= &-\lambda_f^T(t) (f_F)_x(x(t),u(t)) - \lambda_f^T(t) (h_x^T(x(t))z(t))_x \nonumber \\
& &  + \lambda_h(t) h_x(x(t)) \label{eqAdjEqeqAsFirstCaseAdjDaeDiff}\\
0 &= & \lambda_f^T(t) h_x^T(x(t)),\ t\in [t_t^+,t_f]. \label{eqAdjEqsFirstCaseAdjDaeAlg}
\end{eqnarray}
(\ref{eqAdjEqsFirstCaseAdjOde}) is a system of ODEs with respect to the adjoint variable $ \lambda_f(t) $. (\ref{eqAdjEqeqAsFirstCaseAdjDaeDiff})--(\ref{eqAdjEqsFirstCaseAdjDaeAlg}) form a system of DAEs with respect to the adjoint variables $ \lambda_f(t) $ and $ \lambda_h(t) $. The DAEs are of index 2 which can be shown by considering the differentiation of (\ref{eqAdjEqsFirstCaseAdjDaeAlg}) (see, e.g., \cite{ps2020a})
\begin{eqnarray}
0 & = & \left( \lambda_f^T(t) h_x^T(x(t)) \right)' = (\lambda_f^T)'(t) h_x^T(x(t)) + \lambda_f^T(t) \left( h_x^T(x(t)) \right)'\nonumber  \\
& = & \left( -\lambda_f^T(t) (f_F)_x(x(t),u(t)) - \lambda_f^T(t) (h_x^T(x(t))z(t))_x \right. \nonumber \\
& &  \left.  + \lambda_h(t) h_x(x(t)) \right) h_x^T(x(t)) + \lambda_f^T(t) \left( h_x^T(x(t)) \right)'  \nonumber 
\end{eqnarray}
The algebraic variable $ \lambda_h^T(t) $ is multiplied by $ h_x(x(t))h_x^T(x(t)) $ which is nonzero by the assumption {\bf (H1)}. This means that we can evaluate $ \lambda_h(t) $ by using the above equation which proves that the DAEs are index 2 equations. One can also show that under {\bf (H1)} a solution to the equations (\ref{eqAdjEqeqAsFirstCaseAdjDaeDiff})-(\ref{eqAdjEqsFirstCaseAdjDaeAlg}), and to the equations (\ref{eqAdjEqsFirstCaseAdjOde}), are such that $\lambda_f$ is an absolutely continuous functions. This fact justifies the use of the integration by parts formula in our derivations.

Let us consider the component
\begin{eqnarray}
&{\displaystyle \left( \phi_x(x(t_f)) +\lambda_f^T(t_f) \right)y^{x,d}(t_f).}\nonumber 
\end{eqnarray}
At the end of the time interval the trajectory evolves according to DAEs (\ref{eqAdjEqeqAsFirstCaseAdjDaeDiff})-(\ref{eqAdjEqsFirstCaseAdjDaeAlg}). The trajectory, after variation of controls, stays within the hyperplane defined by $ h(x(t)) = 0 $ which means that the variation $ y^{x,d}(t_f) $ is always orthogonal to $ h_x(x(t_f)) $. Therefore, the component disappears if the following condition is satisfied:
\begin{eqnarray}
&{\displaystyle \phi_x(x(t_f)) +\lambda_f^T(t_f) = \nu h_x(x(t_f)).}\nonumber 
\end{eqnarray}
for some real number $\nu$.
 
The DAEs (\ref{eqAdjEqeqAsFirstCaseAdjDaeDiff})-(\ref{eqAdjEqsFirstCaseAdjDaeAlg}) have to be consistently initialized. To provide the consistent endpoint conditions for the adjoint variables $ \lambda_f,\ \lambda_h $ we follow the approach presented in \cite{clps2000} (see the Hessenberg index--2 DAE system considered therein) and solve the following system of equations at time 
$ t_f $ for the variables $ \lambda_f, \lambda_h, \nu $ 
\begin{eqnarray}
\nu h_x^T(x(t_f)) & = & \phi_x^T(x(t_f)) + \lambda_f \label{terminal} \\
0 & = & h_x(x(t_f))\lambda_f  \nonumber \\
0 & = & \left(h_x(x(t_f))\right)'\lambda_f - h_x(x(t_f)) \left( f_F \right)_x^T(x(t_f),u(t_f)) \lambda_f -\nonumber \\
& & h_x(x(t_f)) \left(h_x^T(x(t_f))z(t_f)\right)_x^T \lambda_f +  \nonumber \\
& & h_x(x(t_f)) h_x^T(x(t_f)) \lambda_h. \label{slidingDAEIndex2HiddenAdjEndpoint}
\end{eqnarray}
The solution to these equations is (due to {\bf (H1)}): $\nu = h_x\phi_x^T/\|h_x\|^2$, $\lambda_f = h_x\phi_x^Th_x^T/$ $\|h_x\|^2 - \phi_x^T$, $\lambda_h = \lambda_f^T\left ((\left ( \left (f_F\right )_x + \left (h_xz\right )_x\right )h_x^T + \left (h_x^T\right )'\right )/\|h_x\|^2$, where all functions are evaluated at $(x(t_f),z(t_f),u(t_f))$. Then we set $\lambda_f (t_f) = \lambda_f$, $\lambda_h (t_f) = \lambda_h$.

At the transition time $ t_t $ the adjoint variable $ \lambda_f $ undergoes a jump. To calculate the value of $ \lambda_f(t_t^-) $ the following system of equations have to be solved for the variables $ \lambda_f^-,\ \pi $
\begin{eqnarray}
\lambda_f^- & = & \lambda_f(t_t^+) - \pi h_x^T(x(t_t)) \nonumber  \\
\left (\lambda_f^-\right )^T f_1(x(t_t^-),u(t_t^-)) & = & \lambda_f^T(t_t^+) f_F(x(t_t^+),u(t_t^+)) + \lambda_f^T(t_t^+)h_x^T(x(t_t^+))z(t_t^+) -\nonumber \\
& &  \lambda_h^T(t_t^+)h(x(t_t^+)). \label{odeAdjTransHamilCont}
\end{eqnarray}
The solution to these equations is (due to {\bf (H2)}: $\pi = ( (\lambda_f^+)^Tf_1(-)-$\\ $(\lambda_f^+)^Tf_F(+)-(\lambda_f^+)^Th_x(+)z^+ + \lambda_h^+ h(+))/h_x(+)f_1(-)$, $\lambda_f^- = \lambda_f^+ - \pi h_x(+)$. where $\lambda_f^+$, etc. is a function evaluated at $t_t^+$ (according to our notation (\ref{Jump2})), function $f_F(+)$, etc. is evaluated at $(x(t_t^+),u(t_t^+))$, function $f_1(-)$, etc. at $(x(t_t^-),u(t_t^-))$. Then we set $ \lambda_f(t_t^-) = \lambda_f^-$.

Once we solve the adjoint equations we obtain the adjoint variables $ \lambda_f $ and $ \lambda_h $. Eventually we are in a position to calculate the first variation of the cost function $ \phi(t_f) $ with respect to a control function variation $ d $ as follows
\begin{eqnarray}
&{\displaystyle d\phi(x(t_f)) = d\Phi(x,z,u,\lambda_f,\lambda_h,\pi) = \left \langle \bar{F}_0(u),d\right \rangle = }\nonumber\\
&{\displaystyle -\int_{t_0}^{t_t^-} \lambda_f^T(t) (f_1)_u(x(t),u(t))d(t) dt - \int_{t_t^+}^{t_f} \lambda_f^T(t) (f_F)_u(x(t),u(t))d(t) dt. }\nonumber  
\end{eqnarray} 

The adjoint equations for other cases than {\it Case 1--3} can be derived in the same way. In {\it Appendix \ref{SecAppendixC}} {\it 3--1} is considered.

\section{The weak maximum principle}
\label{SecWMP}

On the basis of the defined adjoint equations we can formulate the weak maximum principle for the considered problem. Suppose that $\bar{u}$ is the problem solution, then, as shown in {\it Theorem \ref{t1}}, $\bar{u}$ will satisfy ${\bf (NC)}$.

The weak maximum principle for the problem {\bf (P)} can assume a quite complicated form depending on the number of switching points triggered by the optimal control $\bar{u}$. In order to exemplify the possible conditions stated by the weak maximum principle consider the case when optimal control $\bar{u}$ forces the system behavior denoted as {\it Case 1--3)} in Section \ref{secAdjoint}.  We call the necessary optimality conditions for that case ${\bf (NC^{13})}$. When writing these conditions we take into account that  (\ref{Jump4}) holds and the fact that $x$ is continuous.

\vspace{5mm}

\noindent ${\bf (NC^{13})}$: There exist: nonnegative numbers $\alpha_j^2$, $j\in I$, numbers $\alpha_i^1$, $i\in E$ such that $\sum_{i\in E}\left | \alpha_i^1\right | +\sum_{j\in I} \alpha_j^2 \neq 0$; piecewise differentiable function $\lambda_f$; piecewise continuous function $\lambda_h$, such that the following hold\\
{\it (i) terminal conditions}---$\lambda_f(t_f)$, $\lambda_h(t_f)$ and $\nu$ are evaluated from equations below
\begin{eqnarray}
-\lambda_f(t_f) & = & \phi_x^T(\bar{x}(t_f)) - \nu h_x^T(\bar{x}(t_f)) +  \nonumber \\
& & \sum_{i\in E} \alpha_i^1\left (g_i^1\right )_x^T (\bar{x}(t_f))+ \sum_{j\in I} \alpha_j^2\left (g_j^2\right )_x^T (\bar{x}(t_f)) \nonumber \\
0 & = & h_x(\bar{x}(t_f))\lambda_f(t_f) \nonumber \\
0 & = & \left(h_x(\bar{x}(t_f))\right)'\lambda_f(t_f) -\nonumber \\
& & 
h_x(\bar{x}(t_f)) \left( f_F \right)_x^T(\bar{x}(t_f),\bar{u}(t_f))\lambda_f(t_f) - \nonumber \\
& &  h_x(\bar{x}(t_f)) \left(h_x^T(\bar{x}(t_f))\bar{z}(t_f)\right)_x^T\lambda_f(t_f)  +  \nonumber \\
& & h_x(\bar{x}(t_f)) h_x^T(\bar{x}(t_f)) \lambda_h (t_f); \nonumber
\end{eqnarray}
{\it (ii) adjoint equations}\\
for $ t\in [t_t,t_f) $
\begin{eqnarray}
\lambda_f' &=& - \left( f_F \right)_x^T(\bar{x},\bar{u})\lambda_f - \left(h_x^T(\bar{x})\bar{z}\right)_x^T\lambda_f + h_x^T(\bar{x}) \lambda_h \nonumber \\ 
0 &=& h_x(\bar{x})\lambda_f \nonumber 
\end{eqnarray}
\noindent for $ t\in [t_0,t_t) $
\begin{eqnarray}
\lambda_f' = - \left( f_1 \right)_x^T(\bar{x},\bar{u})\lambda_f; \nonumber
\end{eqnarray}
{\it (iii) jump conditions}---$\lambda_f(t_t^-)$ and $\pi$ are evaluated from the equations below
\begin{eqnarray}
\lambda_f(t_t^-) & = &\lambda_f(t_t^+) - \pi h_x^T(\bar{x}(t_t))  \nonumber \\
\lambda_f^T(t_t^-) f_1(\bar{x}(t_t),\bar{u}(t_t)) & = & \lambda_f^T(t_t^+) f_F(\bar{x}(t_t),\bar{u}(t_t^+)) + \nonumber \\
& & \lambda_f^T(t_t^+)h_x^T(\bar{x}(t_t))\bar{z}(t_t^+) -  \nonumber \\
& & \lambda_h^T(t_t^+)h(\bar{x}(t_t));
\nonumber
\end{eqnarray}
{\it (iv) the weak maximum principle}\\
\begin{eqnarray}
    &{\displaystyle H^{13}(\bar{x},\bar{u},\lambda_f,t_t,u)\leq H^{13}(\bar{x},\bar{u},\lambda_f,t_t,\bar{u}),\ \forall u\in {\mathcal U}}\nonumber 
\end{eqnarray}
where $H^{13}(\bar{x},\bar{u},\lambda_f,t_t,u) = \int_{t_0}^{t_t^-}\lambda_f^T(t) \left (f_1\right )_u (\bar{x}(t),\bar{u}(t))u(t)dt +$\\ $\int_{t_t^+}^{t_f}\lambda_f^T(t) \left (f_F\right )_u (\bar{x}(t),\bar{u}(t))u(t)dt$.\\
{\it (v) complementarity conditions}
\begin{eqnarray}
\alpha_j^2 = 0,\ \ {\rm if}\ \ j\not\in I_{0,\bar{u}}.
\nonumber 
\end{eqnarray}
{\it (vi) feasibility} 
\begin{eqnarray}
&{\displaystyle 
g_i^1 (\bar{x}(t_f)) = 0,\  i\in E,\  g_j^2 (\bar{x}(t_f))\leq 0,\ j\in I.}\nonumber 
\end{eqnarray}

For the {\it Case 1--3)} we can prove the theorem which is the restatement of part {\it iii)} of {\it Theorem \ref{t1}}.

\newtheorem{Th3}{Theorem}[section]
\begin{Th3}
\label{Th3}
Assume that data for {\rm {\bf (P)}} satisfy hypotheses ${\bf (H1)}$, ${\bf (H2)}$, ${\bf (H3)}$, ${\bf (H4)}$, ${\bf (CC)}$ and ${\bf (CQ)}$. Let $\{u_k\}$ be a sequence of controls generated by {\it Algorithm}. Then 
	\begin{enumerate}
		\item [iii)] if $\bar{u}$ is the limit point of $\{u_k\}$, $\bar{x}$ its state  corresponding trajectory and suppose that for $\bar{u}$ {\it Case 1--3)} applies, then the pair $(\bar{x},\bar{u})$ satisfies the conditions ${\bf (NC^{13})}$.
	\end{enumerate}
\end{Th3}
\begin{proof}
As shown in the proof of Theorem 5.1 ({\it Stage 2}--Dualization) in \cite{pv98} the conclusions of the proof of {\it Theorem \ref{t1}} can be expressed by
\begin{eqnarray}
&{\displaystyle
	\min_{d\in {\mathcal D}_{\bar{u}}} \max_{\gamma\in {\mathcal K}} \Psi(d,\gamma) = 0}\nonumber
\end{eqnarray}
where
\begin{eqnarray}
&{\displaystyle
	{\mathcal K} = \left \{ \gamma = \left (\alpha_0,\{\alpha^1_i\}_{i\in E},
	\{\alpha^2_j\}_{j\in I}\right )\in \mathbb{R}^{1+|E|+|I|}:\   
	\alpha_0\geq 0,\ \alpha^2_j\geq 0,\ j\in I, \right . } \nonumber \\
&{\displaystyle 
	\left . \ \alpha_0 + \sum_{i\in E} \left |\alpha^1_i\right | + 
	\sum_{j\in I}\alpha^2_j = 	1,\ 
	\alpha^2_j = 0\ \ {\rm if}
	\ j\not\in I_{0,\bar{u}} \right \}}\nonumber
\end{eqnarray}
and
\begin{eqnarray}
&{\displaystyle
	\Psi(d,\gamma) := \alpha_0\left \langle \nabla \bar{F}_0(\bar{u}),d\right \rangle + c \left (
	\sum_{i\in E}\alpha^1_i\left \langle \nabla \bar{g}^1_i(\bar{u}),d\right \rangle
	+ \sum_{j\in I_{0,\bar{u}}} \alpha^2_j\left\langle \nabla \bar{g}^2_j(\bar{u}),
	d\right\rangle \right ).}\nonumber
\end{eqnarray}

$\Psi(\cdot,\gamma)$ is a linear function on ${\mathcal C}^N_m[T]$ of which
${\mathcal D}_{\bar{u}}$ is a convex subset.
$\Psi(d,\cdot)$ is a bounded linear map and ${\mathcal K}$ is a compact convex set with
respect to the product topology of $\mathbb{R}^{1+|E|+|I|}$. It follows from the minimax theorem 
(\cite{ae}) 
that there exists some nonzero $\bar{\gamma}\in {\mathcal K}$
such that
\begin{eqnarray}
&{\displaystyle
	\min_{d\in {\mathcal D}_{\bar{u}}}\max_{\gamma\in {\mathcal K}} \Psi(d,\gamma) =
	\min_{d\in {\mathcal D}_{\bar{u}}} \Psi(d,\bar{\gamma}) = 0,}\label{l530}
\end{eqnarray}
with $\bar{\gamma} = (\bar{\alpha}_0,\{\bar{\alpha}^1_i\}_{i\in E},\{\bar{\alpha}^2_j\}_{j\in I})$.

Since the constraint qualification ${\bf (CQ)}$ holds we can easily show that $\bar{\alpha}_0\neq 0$.

The adjoint equations have been derived for the functional $\bar{F}_0(u)$, however similar analysis could be carried out for the functional
\begin{eqnarray}
&{\displaystyle H(u) = \bar{F}_0(u) + c\left ( \sum_{i\in E} \bar{\alpha}_i^1 \bar{g}_i^1 (u) + \sum_{j\in I}\bar{\alpha}_j^2\bar{g}_j^2(\bar{u}) \right ), }\nonumber
\end{eqnarray}	
and we arrive at the optimality conditions by taking $\alpha_i^1 = c\bar{\alpha}_i^1$, $i\in E$, $\alpha_j^2 = c\bar{\alpha}_j^2$, $j\in I$ (notice that $c>0$).
\end{proof}

We do not prove that ${\bf (NC^{13})}$ are necessary optimality conditions for the problem ${\bf (P)}$ when {\it Case 1--3} occurs. However, using arguments similar to these applied in the proofs of {\it Theorem \ref{t1}} and {\it Theorem \ref{Th3}}, techniques of deriving optimality conditions for general optimization problems (\cite{N92},\cite{G72}), by taking into account {\it Lemma \ref{l2}}, we can prove that this is the case. Due to the lack of space the proof of that is not provided in the paper.

The other cases than {\it Case 1--3)} can be analyzed, with respect to necessary optimality conditions, in a similar way. We do not enumerate these cases because they are numerous, e.g.: the system behaves first according to function $f_2$, then, until the end it is in the sliding mode; the system starts in the sliding mode before moving to the mode with dynamics described by $f_1$ or $f_2$; there are several discrete modes the system assumes on the interval $[t_0,t_f]$. In particular, for {\it Case 3--1)} the necessary optimality conditions will shape as stated in {\it  Appendix \ref{SecAppendixD}}.

%% file: sec_withoutSMSC_P2.tex
\section{Hybrid systems without sliding motion} \label{secwSM}
The process of finding optimal solutions to control problems with hybrid systems significantly simplifies when trajectories generated during the process do contain sections with sliding modes. In this case we assume that our state trajectory on the entire horizon $[t_0,t_f]$ will consist of trajectories of type of type $x$ only, i.e. we have 
\begin{eqnarray}
&{\displaystyle [t_0,t_f] = \cup_{i\in I_1}A^1_i\cup_{i\in I_2}A^2_i. } \nonumber 
\end{eqnarray}

For the simplicity of presentation consider the case (we call it {\it Case 1--2)}) when the system on time subinterval $[t_0,t_t)$ evolves according to equations
\begin{eqnarray}
&{\displaystyle x' = f_1(x,u),}\label{wSM1}
\end{eqnarray} 
and then, on the subinterval $[t_t,t_f]$, according to 
\begin{eqnarray}
&{\displaystyle x' = f_2(x,u),}\label{wSM2}
\end{eqnarray}
Due to our previous remarks we do not have to restrict the class of admissible control functions to piecewise smooth functions, and we can can take as ${\mathcal U}$ the set
\begin{eqnarray}
&{\displaystyle {\mathcal U} = \left \{ u\in {\mathcal L}_m^2[T]:\ u(t) \in U,\ {\rm a.e.\ on}\ T\right \},} \label{wSM3}
\end{eqnarray}
For this set of admissible controls, and due to our assumption that sliding mode does not occur, the relations (\ref{Lin1})--(\ref{Lin4}) hold for $x$, $x^d$ and $y^{x,d}$ on the subintervals $[t_0,t_t]$ and $[t_t,t_f]$.  

The adjoint equations for the equations (\ref{wSM1})--(\ref{wSM2}) can be derived in a way described in Section \ref{secAdjoint}. To this end we define the functional

\begin{eqnarray}
&{\displaystyle \Phi(x,u,\lambda_f,\pi) = \phi(x(t_f)) + \pi h(x(t_t))  +}\nonumber \\
&{\displaystyle \int_{t_0}^{t_t^-} \lambda_f^T(t) \left( x'(t) - f_1(x(t),u(t)) \right) dt + \int_{t_t^+}^{t_f}\lambda_f^T(t) \left( x'(t) - f_2(x(t),u(t)) \right) dt,} \nonumber
\end{eqnarray} 
and after applying arguments as in Section \ref{secAdjoint}, we arrive at the equations for adjoint variable $\lambda_f$:
\begin{equation}
(\lambda_f^T)'(t) = - \lambda_f^T(t) (f_1)_x(x(t),u(t)),\ t\in[t_0,t_t]  
\label{eqAdjf1}
\end{equation}
and 
\begin{eqnarray}
&{\displaystyle (\lambda_f^T)'(t) = \lambda_f^T(t) (f_2)_x(x(t),u(t)),\ t\in (t_t,t_f].} \label{eqAdjf2}
\end{eqnarray}

The transversality condition for the equations (\ref{eqAdjf1})--(\ref{eqAdjf2}) is
\begin{eqnarray}
&{\displaystyle \phi_x^T(x(t_f)) + \lambda_f(t_f) = 0} \label{eqAdjterminal} 
\end{eqnarray}
and the jump conditions at the switching time $t_t$ are 
\begin{eqnarray}
\lambda_f(t_t^-) & = & \lambda_f(t_t^+) - \pi h_x^T(x(t_t)) \nonumber  \\
\lambda_f^T(t_t^-) f_1(x(t_t^-),u(t_t^-)) & = & \lambda_f^T(t_t^+) f_2(x(t_t^+),u(t_t^+)). \label{eqAdjswitching}
\end{eqnarray}
Then the first variation of the cost function $ \phi(t_f) $ with respect to a control function variation $ d $ is follows
\begin{eqnarray}
&{\displaystyle d\phi(x(t_f)) = d\Phi(x,z,u,\lambda_f,\pi) = }\nonumber\\
&{\displaystyle -\int_{t_0}^{t_t^-} \lambda_f^T(t) (f_1)_u(x(t),u(t))d(t) dt  - \int_{t_t^+}^{t_f} \lambda_f^T(t) (f_2)_u(x(t),u(t))d(t) dt. }\nonumber  
\end{eqnarray} 

Having adjoint equations we can formulate necessary optimality conditions for the case---${\bf (NC^{12})}$. 
\vspace{3mm}

\noindent ${\bf (NC^{12})}$: There exist: nonnegative numbers $\alpha_j^2$, $j\in I$, numbers $\alpha_i^1$, $i\in E$ such that $\sum_{i\in E}\left | \alpha_i^1\right | +\sum_{j\in I} \alpha_j^2 \neq 0$; absolutely continuous function $\lambda_f$, such that the following hold:\\
{\it (i) terminal conditions}
\begin{eqnarray}
\lambda_f(t_f) & = & -\phi_x^T(\bar{x}(t_f))  \nonumber 
\end{eqnarray}
{\it (ii) adjoint equations}\\
for $ t\in (t_t,t_f) $
\begin{eqnarray}
\lambda_f' &=& - \left( f_2 \right)_x^T(\bar{x},\bar{u})\lambda_f\nonumber 
\end{eqnarray}
\noindent for $ t\in [t_0,t_t] $
\begin{eqnarray}
\lambda_f' = - \left( f_1 \right)_x^T(\bar{x},\bar{u})\lambda_f; \nonumber
\end{eqnarray}
{\it (iii) jump conditions}---$\lambda_f(t_t^-)$ and $\pi$ are the solutions to the algebraic equations stated below
\begin{eqnarray}
\lambda_f(t_t^-) & = &\lambda_f(t_t^+) - \pi h_x^T(\bar{x}(t_t^-))  \nonumber \\
\lambda_f^T(t_t^-) f_1(\bar{x}(t_t^-),\bar{u}(t_t^-)) & = & \lambda_f^T(t_t^+) f_2(\bar{x}(t_t^+),\bar{u}(t_t^+)) 
\nonumber
\end{eqnarray}
{\it (iv) the weak maximum principle}\\
\noindent for $ t\in (t_t,t_f) $
\begin{eqnarray}
\lambda_f^T(t) \left (f_2\right )_u (\bar{x}(t),\bar{u}(t))u \leq 
\lambda_f^T(t) \left (f_2\right )_u (\bar{x}(t),\bar{u}(t))\bar{u}(t)\label{withoutSDMCa} 
\end{eqnarray}
\noindent for $ t\in [t_0,t_t] $
\begin{eqnarray}
\lambda_f^T(t) \left (f_1\right )_u (\bar{x}(t),\bar{u}(t))u \leq 
\lambda_f^T(t) \left (f_1\right )_u (\bar{x}(t),\bar{u}(t))\bar{u}(t)\label{withoutSDMCb} 
\end{eqnarray}
for all $u\in U$;\\
{\it (v) complementarity conditions}
\begin{eqnarray}
\alpha_j^2 = 0,\ \ {\rm if}\ \ j\not\in I_{0,\bar{u}}.
\nonumber 
\end{eqnarray}
{\it (vi) feasibility} 
\begin{eqnarray}
&{\displaystyle 
g_i^1 (\bar{x}(t_f)) = 0,\  i\in E,\  g_j^2 (\bar{x}(t_f))\leq 0,\ j\in I.}\nonumber 
\end{eqnarray}

Since we elaborate on the case of optimal control problems with hybrid systems which do not involve sliding modes we have to introduce the assumption which excludes the possibility of system going into the sliding mode. The assumption concerns switching points $t_t$ and is postulated instead of ${\bf (H2)}$.

\begin{description}
	\item[${\bf (H2')}$]Function $h(\cdot)$ is differentiable and there exist $0<L_1<+\infty$ and $0<L_2 <+\infty$ such that
	\begin{eqnarray}
	&{\displaystyle 
		\left \| h_{x}(\hat{x}) - h_{x}(x)\right \| \leq  L_1\|\hat{x} - x \|,} \label{H3a12} \\
	&{\displaystyle 
		\left | h_x(\hat{x})f_{i}(\hat{x},\hat{u}) -  h_x(x)f_{i}(x,u)\right | \leq  L_2\|(\hat{x},\hat{u}) - (x,u) \|,} \label{H3b12}	
	\end{eqnarray}
	for all  $(x,u)$, $(\hat{x},\hat{u})$ in $\mathbb{R}^n \times U$.

For any $u\in {\mathcal U}$ and any switching point $t_t$ the following limits exist
		\begin{eqnarray}
		&{\displaystyle \lim_{t\rightarrow t_t, t< t_t} u(t),\ \lim_{t\rightarrow t_t, t> t_t} u(t)}
		\end{eqnarray}
		(and are denoted by $u(t_t^-)$ and $u(t_t^+)$ respectively).
		
	There exists $\varepsilon > 0$ and $0<L_3 <+\infty$ such that for all switching points $t_t$ (each switching time corresponds to some $u\in {\mathcal U}$) and for all their perturbations $t_t^d$ triggered by perturbations $d$ such that $u+d\in {\mathcal U}$, $\|d\|_{{\mathcal L}^2}\leq \varepsilon$, and for all $\theta\in [0,1]$ we have
	\begin{eqnarray}
	&{\displaystyle 
		h_x(x(\tau^{\theta ,d}(t_t)) ) f_{i}(x(\tau^{\theta ,d}(t_t)),u(\tau^{\theta ,d}(t_t) ))\geq L_3,} \label{H3c12}	
	\end{eqnarray}
	or
	\begin{eqnarray}
&{\displaystyle 
	h_x(x(\tau^{\theta ,d}(t_t)) ) f_{i}(x(\tau^{\theta ,d}(t_t)),u(\tau^{\theta ,d}(t_t) )) \leq -L_3,} \label{H3cc12}	
\end{eqnarray}	
for $i=1,2$, 	
	where
	\begin{eqnarray}
	\tau^{\theta ,d}(t_t) & = & t_t + \theta (t^d_t - t_t). 
	\label{Midt12} 
	\end{eqnarray}

\end{description}

The hypothesis {\bf (H2')} which replaces {\bf (H2)} does not request the condition (\ref{Dprime}) for perturbations $d$. This is the consequence of having the switching function $\eta=h$ which does not depend on controls. That simplification spreads through the definitions which concern perturbations $d$, for example, the definitions of ${\mathcal U}$, ${\mathcal D}_u$ and ${\mathcal D}$, should be updated in order to work with the broader sets of perturbations. In particular, in the hypothesis ${\bf (CC)}$ the requirement (\ref{Dprime_equiv}) is not needed.  {\it Theorem \ref{Th4}} can be proved along the lines of the proof of {\it Theorem \ref{t1}}, {\it Theorem \ref{Th3}} and with sets ${\mathcal D}_u$, ${\mathcal D}$ modified in this way. 

\newtheorem{Th4}{Theorem}[section]
\begin{Th4}
\label{Th4}
Assume that data for {\rm {\bf (P)}} satisfy hypotheses ${\bf (H1)}$ (part i)),\\ ${\bf (H2')}$, ${\bf (H3)}$, ${\bf (H4)}$, ${\bf (CC)}$ and ${\bf (CQ)}$. Let $\{u_k\}$ be a sequence of controls generated by {\it Algorithm} and let $\{c_k\}$ be a sequence of the corresponding penalty parameters. Then 
	\begin{enumerate}
		\item [i)] $\{c_k\}$ is a bounded sequence
		\item [ii)] $\lim_{k\rightarrow\infty}\sigma_{c_k}(u_k)=0$, $\lim_{k\rightarrow\infty} M(u_k) = 0$,
		\item [iii)] if $\bar{u}$ is the limit point of $\{u_k\}$, $\bar{x}$ its state  corresponding trajectory and suppose that for $\bar{u}$ Case 1--2) applies, then the pair $(\bar{x},\bar{u})$ satisfies the conditions ${\bf (NC^{12})}$.
	\end{enumerate}
\end{Th4}
\begin{proof}
Following analysis applied in the proofs of {\it Theorem \ref{Th3}} and {\it Theorem \ref{t1p}} we can show that if $\bar{u}$ is a local minimizer to the problem ${\bf (P)}$ and the assumptions of the theorem hold then
\begin{eqnarray}
&{\displaystyle H^{12}(\bar{x},\bar{u},\lambda_f,t_t,u)\leq H^{12}(\bar{x},\bar{u},\lambda_f,t_t,\bar{u}),\ \forall u\in {\mathcal U}}\nonumber 
\end{eqnarray}
where 
\begin{eqnarray}
&{\displaystyle H^{12}(\bar{x},\bar{u},\lambda_f,t_t,u) = \int_{t_0}^{t_t^-}\lambda_f^T(t) \left (f_1\right )_u (\bar{x}(t),\bar{u}(t))u(t)dt}\nonumber \\
&{\displaystyle +\int_{t_t^+}^{t_f}\lambda_f^T(t) \left (f_2\right )_u (\bar{x}(t),\bar{u}(t))u(t)dt.}
\end{eqnarray}

Referring to the mean value theorem for integrals and taking into account the definition of the admissible set of controls (\ref{wSM3}) (the set $U$ is bounded) we come to the conditions (\ref{withoutSDMCa})--(\ref{withoutSDMCb})
\end{proof}

%% file: sec_conclusionsSC_P2.tex
\section{Conclusions}
\label{SecConclusions}

The paper presents the computational approach to hybrid optimal control problems with sliding modes. It seems to be the first method for optimal control problems with hybrid systems which can exhibit sliding modes. We show that under  the assumption that controls are represented by piecewise smooth functions the presented method is globally convergent in the sense that every accumulation point of a sequence generated by the method satisfies necessary optimality conditions in the form of the weak maximum principle. In the accompanying papers (\cite{ps19a},\cite{ps19b}) we present the implementable version of our computational method by relaxing the condition that continuous state trajectories are available. In these papers we assume that system equations are integrated by an implicit Runge--Kutta method and we show that the global convergence of the implementable version of our algorithm is preserved provided that the integration procedure step sizes converge to zero.

%% file: sec_appendixSC_P2_ver2.tex
\appendix
\label{SecAppendix}
\section{Proof of {\it Proposition \ref{p1}}}
\label{AppendixP}
\begin{proof} {\it (i)} Fix $u\in {\mathcal U}$,  we must find $\hat{c} > 0$ such that, if $c > \hat{c}$ then
	$t_c(u)\leq 0$. If $M(u)=0$ then of course
	$t_c(u)\leq 0$ for any $c >0$. If $M(u) > 0$ then according to {\it Lemma \ref{l2}}
	there exists
	$\hat{d}\in {\cal U}-u$ such that, if we set $\varepsilon=K_1M(u) > 0$ 
	with $K_1$ as in {\it Lemma \ref{l2}}, then
	\begin{eqnarray}
	&{\displaystyle 
		\max_{i\in E} \left [\left |\bar{g}^1_i(u) + \left \langle \nabla \bar{g}^1_i(u),\hat{d} \right \rangle \right | -M(u) \right ] \leq 
		\max_{i\in E} \left [\left | \bar{g}^1_i(u) + \left \langle \nabla \bar{g}^1_i(u),\hat{d} \right \rangle \right | - \right .}\nonumber \\
  &{\displaystyle \left . \left | \bar{g}_i^1(u) \right | \right ] \leq  -\varepsilon,} \nonumber \\
	&{\displaystyle 
		\max_{j\in I} \left [\bar{g}^2_j(u) + \left \langle \nabla \bar{g}^2_j(u),\hat{d} \right \rangle  -M(u) \right ] \leq
		\max_{j\in I} \left [\bar{g}^2_j(u) + \left \langle \nabla \bar{g}^2_j(u),\hat{d} \right \rangle  - \right . }\nonumber \\
  &{\displaystyle \left . \bar{g}^2_i(u) \right ]  \leq  -\varepsilon } \nonumber 
	\end{eqnarray}
	and 
	\begin{eqnarray}
	&{\displaystyle
		\theta (u) < -\varepsilon.}\nonumber
	\end{eqnarray}
	Here
	\begin{eqnarray}
	&{\displaystyle
		\theta(u) = \max\left [ \max_{i\in E} \left |\bar{g}^1_i(u) + 
		\left \langle \nabla \bar{g}^1_i(u),\hat{d}\right \rangle \right |
		-M(u), \right .}
	\nonumber \\
	&{\displaystyle
		\left . 
		\max_{j\in I}\left [ \bar{g}^2_j(u) + \left \langle \nabla \bar{g}^2_j(u),
		\hat{d}\right \rangle - M(u) \right ]\right ].}
	\nonumber 
	\end{eqnarray}
	Because $\sigma_{c}(u) = \left \langle \nabla \bar{F}_0(u),\bar{d}\right \rangle + 
	c\left [\bar{\beta}-M(u)\right ]$,
	from the definition of $t_c$ and (\ref{Lin1})--(\ref{Lin3}),
	we get  
	\begin{eqnarray}
	&{\displaystyle
		t_c(u) \leq W + c \theta(u) + M(u)/c,}\nonumber
	\end{eqnarray}
	where $W=\max\ [0,\langle \nabla \bar{F}_0(u),\hat{d}\rangle]$.
	It follows that $t_c(u)\leq 0$ for any $c > \hat{c}$ where
	\begin{eqnarray}
	&{\displaystyle
		\hat{c} := \max \left [1,  \frac{W + M(u)}{-\theta(u)}\right ]. }\nonumber
	\end{eqnarray}
	
	{\it (ii)} Take $u\in {\mathcal U}$ and $c > 0$ such that $\sigma_c(u) < 0$. Let
	$(d,\beta)$ be the solution to ${\bf P_c(u)}$. 
	Since $\sigma_c(u)\neq 0$, it follows that $d\neq 0$.    
	
	We deduce from the differentiability properties of $\phi$, $g^1_i$, $g^2_j$, and (\ref{Lin4})
	that there exists $o:[0,\infty)\rightarrow [0,\infty)$ such that $s^{-1} o(s)\rightarrow
	0$ as $s\downarrow 0$ and the following inequality is valid for any $\alpha\in
	[0,1]$:
	\begin{eqnarray}
	&{\displaystyle 
		\bar{F}_c(u+\alpha d) - \bar{F}_c(u) \leq  
		\alpha \left \langle \nabla \bar{F}_0(u),d\right \rangle +
		c \max \left [\max_{i\in E} \left | \bar{g}^1_i(u) + \right . \right .}\nonumber \\
  &{\displaystyle \left . \left . \alpha \left \langle \nabla \bar{g}^1_i(u),
		d\right \rangle \right |, \right . 
		\left . \max_{j\in I} \left ( \bar{g}^2_j(u) + \alpha \left \langle \nabla \bar{g}^2_j(u),
		d\right \rangle \right ) \right ] - M(u) + o(\alpha).} \label{l59}
	\end{eqnarray}
	
	By convexity of the functions $e\rightarrow \max_{i\in E} | \bar{g}^1_i(u) +
	\langle \nabla \bar{g}^1_i(u),e\rangle |$, $e\rightarrow \max_{j\in I} ( \bar{g}^2_j(u)$ $+\langle \nabla \bar{g}^2_j(u),e\rangle )$: 
	\begin{eqnarray}
	&{\displaystyle
		\max_{i\in E} \left | \bar{g}^1_i(u) + \alpha \left \langle \nabla \bar{g}^1_i(u),
		d\right \rangle \right | 
		- M(u)
		\leq 
		\alpha \left [ \max_{i\in E} \left | \bar{g}^1_i(u) +
		\left \langle \nabla \bar{g}^1_i(u),d\right \rangle \right |\right . }\nonumber \\
		&{\displaystyle \left .- M(u) 
		\right ],}\nonumber \\
	&{\displaystyle
		\max_{j\in I} \left ( \bar{g}^2_j(u) + \alpha \left \langle \nabla \bar{g}^2_j(u),
		d\right \rangle \right )
		- M(u)
		\leq 
		\alpha \left [ \max_{j\in I} \left ( \bar{g}^2_j(u) +
		\left \langle \nabla \bar{g}^2_j(u),d\right \rangle \right ) \right . }\nonumber \\
		&{\displaystyle \left . - M(u)
		\right ].}\nonumber
	\end{eqnarray}
	From inequality (\ref{l59}) then
	\begin{eqnarray}
	&{\displaystyle 
		\bar{F}_c(u+\alpha d) - \bar{F}_c(u) \leq \alpha 
		\left \{ \left \langle \nabla \bar{F}_0(u),
		d\right \rangle  + c\max [\max_{i\in E} \left |\bar{g}^1_i(u) + \left \langle\nabla \bar{g}^1_i(u),
		d\right \rangle\right |, 
		\right .}
	\nonumber \\
	&{\displaystyle 
		\left . \left . \max_{j\in I} \left (\bar{g}^2_j(u) + \left \langle\nabla \bar{g}^2_j(u),
		d\right \rangle\right )
		\right ] - c M(u) \right \} + o(\alpha) 
		\leq \alpha \sigma_c(u) + o(\alpha).}\nonumber
	\end{eqnarray}
	It follows that
	\begin{eqnarray}
	&{\displaystyle
		\bar{F}_c(u + \alpha d) - \bar{F}_c(u) \leq \alpha \gamma \sigma_c(u)\ \ \forall
		\alpha\in [0,\alpha_1], }\label{l515}
	\end{eqnarray}
	where $\alpha_1 > 0$ is such that $o(\beta)\leq \beta(\gamma-1)\sigma_c(u)$ for all
	$\beta\in [0,\alpha_1]$.
\end{proof}

\section{Proof of {\it Theorem \ref{t1}}}
\label{AppendixT}
\begin{proof} 
	{\it (i)} Let $\{u_k\}$ be the sequence generated
	by {\it Algorithm} and let $\{c_k\}$ be the corresponding penalty parameters.
	Let $\{k_l\}$ be the sequence of index
	values at which the penalty parameter increases. By extracting a further subsequence (we do not
	relabel) we can arrange that the sequence $\{u_{k_l}\}$ has
	a limit point $\bar{u}\in {\mathcal U}$ because ${\mathcal U}$ is compact.
	We shall find a number
	$\hat{c}<\infty$ such that for sufficiently large $k_l$ `$c_{k_l} > \hat{c}$' implies
	`$\sigma_{c_{k_l}}(u_{k_l})\leq - M(u_{k_l})|/c_{k_l}$'.
	This contradicts our assumption that the penalty parameter increases along the
	subsequence. So we may conclude that $\{c_k\}$ is bounded.
	
	Fix $k_l$ such that $u_{k_l}\in {\mathcal O}(\bar{u})$, where ${\mathcal O}(\bar{u})$ is the
	neighborhood of $\bar{u}$ as specified in {\it Lemma \ref{l2}}.
	From the minimizing property of $\sigma_{c_{k_l}}(u_{k_l})$ we deduce
	\begin{eqnarray}
	\sigma_{c_{k_l}}(u_{k_l}) & \leq &  1/2 \|v_{k_l}-u_{k_l}\|^2_{{\mathcal L}^2} + 
	\left \langle \nabla \bar{F}_0(u_{k_l}),v_{k_l}-u_{k_l}\right \rangle +
	\nonumber \\
	& & c_{k_l} \max\left [ \max_{i\in E} \left | \bar{g}^1_i(u_{k_l}) + 
	\left \langle \nabla \bar{g}^1_i(u_{k_l}),
	v_{k_l} - u_{k_l} \right \rangle \right | - M(u_{k_l}),\right .
	\nonumber \\
	& & \left . 
	\max_{j\in I}\left [\bar{g}^2_j(u_{k_l}) + 
	\left \langle \nabla \bar{g}^2_j(u_{k_l}),v_{k_l}-u_{k_l}\right \rangle \right ] - M(u_{k_l})\right ] 
	\label{equal}
	\end{eqnarray}
	for a control function $v_{k_l}$ satisfying conditions 
	(\ref{l532})--(\ref{l535}) of {\it Lemma \ref{l2}}   
	in which $v_{k_l}$, $u_{k_l}$ replace $v$, $u$ respectively. 
	It follows
	\begin{eqnarray}
	\sigma_{c_{k_l}} (u_{k_l}) & \leq & 1/2\|v_{k_l}-u_{k_l}\|^2_{{\mathcal L}^2} + 
	\left |\left \langle \nabla \bar{F}_0(u_{k_l}),v_{k_l}-u_{k_l}\right \rangle \right |
	- c_{k_l} K_1 M(u_{k_l}).
	\nonumber
	\end{eqnarray}
	Since the control constraint $U$ is bounded and in view of (\ref{Lin1})--(\ref{Lin3})
	there exists $r>0$ (independent of $k_l$) such that
	\begin{eqnarray}
	&{\displaystyle
		\left | \left \langle \nabla
		\bar{F}_0(u_{k_l}),v_{k_l}-u_{k_l}\right \rangle \right |
		\leq r \| v_{k_l} - u_{k_l}\|_{{\mathcal L}^\infty}.}
	\nonumber
	\end{eqnarray}
	Hence
	\begin{eqnarray}
	&{\displaystyle
		\sigma_{c_{k_l}}(u_{k_l}) \leq - (c_{k_l}K_1 - r K_2) M(u_{k_l}).}
	\nonumber
	\end{eqnarray}
	We conclude that
	\begin{eqnarray}
	&{\displaystyle
		\sigma_{c_{k_l}}(u_{k_l}) \leq - M(u_{k_l})/c_{k_l} }
	\nonumber
	\end{eqnarray}
	if $c_{k_l} \geq \hat{c}$, where $\hat{c}$ is the smallest positive number $c$ which satisfies the equation: $c^2 K_1 - c (rK_2) \geq 1$.
	
	\vspace{2mm}
	{\it (ii)} and {\it (iii)} 
	Let $\{u_k\}$ be an infinite sequence generated by {\it Algorithm}.
	We must show that  $\lim_{k\rightarrow\infty}\sigma_{c_k}(u_k) =0$  and,
	if a convergent subsequence of $\{u_k\}$ has a limit point $\bar{u}\in {\mathcal U}$,
	that conditions ${\bf (NC)}$ are satisfied at $\bar{u}$.
	
	Since the $c_k$'s are
	bounded and can increase only by multiplies of $c^0$, we must have $c_k=c$ for all
	$k\geq k_0$, for some $k_0$ and $c>0$. In view of the manner in which $u_k$'s are
	constructed, we have
	\begin{eqnarray}
	&{\displaystyle
		\tilde{F}_c(u_{k+1})-\tilde{F}_c(u_k) \leq
		\gamma \alpha_k \sigma_c(u_k) }\nonumber
	\end{eqnarray}
	for all $k\geq k_0$. This means that, for all $j\geq 1$, $k\geq k_0$
	\begin{eqnarray}
	&{\displaystyle
		\bar{F}_c(u_{k+j}) - \bar{F}_c(u_{k}) \leq \gamma \sum_{i=0}^{j-1}
		\alpha_{k+i} \sigma_c(u_{k+i}). }\label{54b}
	\end{eqnarray}
	Since $\{\bar{F}_c(u_k)\}$ is a bounded sequence and
	$\alpha_k\sigma_c(u_k)$ are nonpositive, we conclude
	\begin{eqnarray}
	&{\displaystyle
		\alpha_k\sigma_c(u_k) \rightarrow 0\ \ {\rm as}\ k\rightarrow \infty.}\label{l520}
	\end{eqnarray}
	
	Since $\sigma_c(u)$ is bounded as $u$ ranges over ${\mathcal U}$ we can arrange
	by a subsequence extraction (we do not relabel) that
	\begin{eqnarray}
	&{\displaystyle
		\sigma_c(u_k) \rightarrow \sigma\ {\rm for\ some}\ \sigma \leq 0.}\nonumber
	\end{eqnarray}
	We claim that $\sigma=0$. To show this, suppose to the contrary that
	$\sigma < 0$. Then by (\ref{l520}) $\alpha_k\rightarrow 0$. 
	
	We must have 
	\begin{eqnarray}
	\bar{F}_c(u_k + \eta^{-1}\alpha_k d_k) - \bar{F}_c (u_k) & > &  \gamma \eta^{-1}\alpha_k
	\sigma_c(u_k)
	\nonumber 
	\end{eqnarray}
	for all $k$ sufficiently  large.                                                
	
	However $\| \alpha_k d_k \|_{{\mathcal L}^\infty}\rightarrow 0$ as
	$k\rightarrow\infty$ (it follows from the presence of term $1/2 \| d\|^2_{{\mathcal L}^2}$ in the objective function of the direction finding subproblem). 
	
	It follows now from (\ref{Lin1})--(\ref{Lin4}) and Propositions 2.1.5--2.1.6 stated in \cite{py99}
	that there exists a function $o:[0,\infty)
	\rightarrow [0,\infty)$ such that $s^{-1}o(s)\rightarrow 0$ as $s\downarrow 0$ and
	\begin{eqnarray}
	\max_{i\in E} \left |\bar{g}^1_i(u_k + \eta^{-1}\alpha_k d_k)\right | & \leq &
	\max_{i\in E}\left | \bar{g}^1_i(u_k) + \left \langle \nabla \bar{g}^1_i(u_k),\eta^{-1}\alpha_k
	d_k\right \rangle \right | + 
	\nonumber \\
	& & o(\eta^{-1}\alpha_k \left \| d_k\right \|_{{\mathcal L}^\infty}),
	\label{l521d} \\
	\max_{j\in I} \bar{g}^2_j(u_k + \eta^{-1}\alpha_k d_k) & \leq &
	\max_{j\in I}\left [ \bar{g}^2_j(u_k) + \left \langle \nabla \bar{g}^2_j(u_k),\eta^{-1}\alpha_k
	d_k\right \rangle \right ] + 
	\nonumber \\
	& & o(\eta^{-1}\alpha_k \left \| d_k\right \|_{{\mathcal L}^\infty}),
	\label{l521dd}
	\end{eqnarray}
	and
	\begin{eqnarray}
	&{\displaystyle 
		\bar{F}_c(u_k+\eta^{-1}\alpha_kd_k) - \bar{F}_c(u_k) \leq 
		\left \langle \nabla \bar{F}_0(u_k),\eta^{-1}\alpha_kd_k\right \rangle +}
	\nonumber \\
	&{\displaystyle 
		c \max \left [\max_{i\in E}
		\left |\bar{g}^1_i(u_k) + \left \langle \nabla \bar{g}^1_i(u_k),\eta^{-1}
		\alpha_kd_k\right \rangle\right | \right .
		-M(u_k),} \nonumber
	\\
	&{\displaystyle 
		\left . \max_{j\in I}
		\bar{g}^2_j(u_k) + \left \langle \nabla \bar{g}^2_j(u_k),\eta^{-1} \alpha_kd_k\right \rangle
		-M(u_k)\right ] + 
		o(\eta^{-1}\alpha_k\left \| d_k\right \|_{{\mathcal L}^\infty}).}
	\label{l521f}
	\end{eqnarray}
	
	By convexity of the functions $e\rightarrow \max_{i\in E} | \bar{g}^1_i(u) +
	\langle \nabla \bar{g}^1_i(u),e\rangle |$, $e\rightarrow \max_{j\in I} ( \bar{g}^2_j(u)$ $+\langle \nabla \bar{g}^2_j(u),e\rangle )$:     
	\begin{eqnarray}
	&{\displaystyle
		\max_{i\in E} \left | \bar{g}^1_i(u_k) + \left \langle \nabla \bar{g}^1_i(u_k),\eta^{-1}\alpha_k
		d_k\right \rangle \right | - 
		M(u_k) }
	\nonumber \\
	&{\displaystyle
		\leq \eta^{-1} \alpha_k \left ( \max_{i\in E}\left |\bar{g}^1_i(u_k) + \left \langle \nabla
		\bar{g}^1_i(u_k),d_k\right \rangle \right | -
		M(u_k) \right ),}
	\label{l521i} \\
	&{\displaystyle
		\max_{j\in I} \bar{g}^2_j(u_k) + \left \langle \nabla \bar{g}^2_j(u_k),\eta^{-1}\alpha_k
		d_k\right \rangle - M(u_k) }
	\nonumber \\
	&{\displaystyle
		\leq \eta^{-1} \alpha_k \left ( \max_{j\in I}\bar{g}^2_j(u_k) + \left \langle \nabla
		\bar{g}^2_j(u_k),d_k\right \rangle -
		M(u_k)\right ),}
	\label{l521ia} 
	\end{eqnarray}
	
	Combining inequalities (\ref{l521d})--(\ref{l521ia}), noting the definition of $\sigma_c(u_k)$
	and the fact that $(d_k,\beta_k)$ solves ${\bf P_{c}(u_k)}$, and dividing across
	the resulting inequality by $\alpha_k$ we arrive at
	\begin{eqnarray}
	&{\displaystyle
		\eta^{-1}\sigma_c(u_k) + \alpha_k^{-1} o(\eta^{-1}\alpha_k \left \| d_k\right \|_{{\mathcal L}^\infty}) \geq \eta^{-1}\gamma \sigma_c(u_k).}\nonumber
	\end{eqnarray}
	We get $\eta^{-1}\sigma \geq \eta^{-1}\gamma \sigma$ in the limit. But
	this implies $\gamma \geq 1$, since $\sigma< 0$ by assumption. From this
	contradiction we conclude the validity of $\sigma =0$. 
	Assertion {\it (ii)} of the theorem
	follows from
	the definition of $t_c$ and part {\it (i)}.
	
	Let $\{u_k\}$ be a convergent subsequence
	with the limit point $\bar{u}\in {\mathcal U}$.
	We must show that conditions ${\bf (NC)}$
	are satisfied at $\bar{u}$.                                                            
	First we establish that $\bar{u}$ is feasible for ${\bf (P)}$ and, for some $c> 0$,
	\begin{eqnarray}
	&{\displaystyle
		0 \leq \left [ \left \langle \nabla \bar{F}_0(\bar{u}),d \right \rangle + c \max \left [\max_{i\in E} \left |\left\langle \nabla \bar{g}^1_i(\bar{u}),d\right\rangle\right|,
		\max_{j\in I} \left\langle \nabla \bar{g}^2_j(\bar{u}),d\right\rangle \right ] \right ]
	} \nonumber
	\end{eqnarray}
	for all $d\in {\mathcal D}_{\bar{u}}$.
	
	Since $\{u_k\}\rightarrow\bar{u}$
	we know that $x_k=x^{u_k}\rightarrow x^{\bar{u}}=\bar{x}$ uniformly.
	Because the penalty parameter is not updated for $k$ sufficiently large we have
	\begin{eqnarray}
	&{\displaystyle
		\sigma_c(u_k) \leq - M(u_k)/c.}\nonumber
	\end{eqnarray}
	But we have shown that $\sigma_c(u_k)\rightarrow 0$ as $k\rightarrow\infty$. It follows
	now from the fact that $u_k\rightarrow\bar{u}\in {\mathcal U}$ 
	that in the limit
	\begin{eqnarray}
	&{\displaystyle
		\max_{i\in E}\left |\bar{g}^1_i(\bar{u})\right | = 0,
		\ \max_{j\in I}\bar{g}^2_j(\bar{u})\leq 0.}\label{limit}
	\end{eqnarray}
	We have established that $\bar{u}$ is feasible for ${\bf (P)}$.
	
	Since $\sigma_c(u_k)\rightarrow 0$ and $M(u_k)\rightarrow 0$ we must also have $d_k\rightarrow 0$. If it does not happen then we would come to the contradiction with (\ref{SigmaPos}). That implies that we also must have $\beta_k\rightarrow 0$ and that means that at $\bar{u}$ (\ref{NC1b})--(\ref{NC1c}) must hold.
	
	Now choose any sequence $\rho_k\downarrow 0$, $\rho_k\leq 1$ for all $k$ such that
	\begin{eqnarray}
	&{\displaystyle
		\rho_k^{-1}\left [ \sigma_c(u_k) - M(u_k)\right ] \rightarrow 0\ \ {\rm as}\ k\rightarrow \infty.}\label{l523}
	\end{eqnarray}
	Choose any $\tilde{d}_k\in {\mathcal D}_{u_k}$. 
	By the convexity of ${\mathcal U}$, $\rho_k \tilde{d}_k\in {\mathcal D}_{u_k}$ for each $k$.
	
	By definition of $\sigma_c(u_k)$ then
	\begin{eqnarray}
	\sigma_c(u_k) & \leq & 1/2 \rho^2_k \|\tilde{d}_k\|^2_{{\mathcal L}^2} + \phi_x(x^{u_k}(t_f)) y^{x_k,\rho_k \tilde{d}_k}(t_f) +
	\nonumber \\
	& & c\max\left [\max_{i\in E}\left | g^1_i(x^{u_k}(t_f)) + (g^1_i)_x (x^{u_k}(t_f))
	y^{x_k,\rho_k \tilde{d}_k}(t_f) \right | - M(u_k),\right .
	\nonumber \\
	& & \left . \max_{j\in I}\left [g^2_j(x^{u_k}(t_f)) + (g^2_j)_x (x^{u_k}(t_f))
	y^{x_k,\rho_k \tilde{d}_k}(t_f)\right ] - M(u_k)\right ].
	\nonumber \\
	\label{l524}
	\end{eqnarray}
	
	Fix $\hat{\varepsilon} > 0$. Since $\rho_k\downarrow 0$ (and consequently
	$y^{x_k,\rho_k\tilde{d}_k}\rightarrow 0$ uniformly),
	we have:
	\begin{eqnarray}
	&{\displaystyle
		\max_{j\in I}\left [ g^2_j(x^{u_k}(t_f)) + (g^2_j)_x (x^{u_k}(t_f))
		y^{x_k,\rho_k\tilde{d}_k}(t_f)\right ]}
	\nonumber \\
	&{\displaystyle
		\leq \max_{I_{\hat{\varepsilon},\bar{u}}} \left [
		\ g^2_j(x^{u_k}(t_f)) + (g^2_j)_x(x^{u_k}(t_f)) y^{x_k,\rho_k\tilde{d}_k}(t_f) \right ], }
	\nonumber
	\end{eqnarray}
	for all $k$ sufficiently large. Here, $I_{\varepsilon,u} = \{i\in I:\ \bar{g}_i^2(u) \geq \max_{j\in I} [\bar{g}_j^2(u)]-\varepsilon \}$. Inserting these inequalities into (\ref{l524}),
	noting that 
	$y^{x_k,\rho_k \tilde{d}_k}=\rho_ky^{x_k,\tilde{d}_k}$, dividing across by
	$\rho_k$ and passing to the limit 
	with the help of (\ref{l523}) and continuity of $\bar{F}_0(\cdot)$, $u\rightarrow
	\langle \nabla\bar{F}_0(u),d\rangle$, etc, we obtain
	\begin{eqnarray}
	&{\displaystyle
		0\leq \phi_x(x^{\bar{u}}_r(t_f)) y^{\bar{x},d}(t_f) +c\max\left [\max_{i\in E}\left |(g^1_i)_x
		(x^{\bar{u}}(t_f)) y^{\bar{x},d}(t_f)\right|,
		\right . }
	\nonumber \\
	&{\displaystyle
		\left .
		\max_{j\in I_{\hat{\varepsilon},\bar{u}}}(g^2_j)_x(x^{\bar{u}}(t_f)) y^{\bar{x},d}(t_f) \right ]. }
	\label{mpwlm}
	\end{eqnarray}
	This inequality is valid for each $\hat{\varepsilon} > 0$ and $d\in {\mathcal D}_{\bar{u}}$.
	
	Again choose arbitrary $d\in {\cal D}_{\bar{u}}$ and 
	take $\varepsilon_k\downarrow 0$. For each $k$ let
	$(g^2_j)_x(x^{\bar{u}}(t_f))\cdot$ $y^{\bar{x},d}(t_f)$ achieves its maximum over
	$I_{\varepsilon_k,\bar{u}}$ at $j=j_k$.
	Then
	\begin{eqnarray}
	&{\displaystyle
		0 \leq  \phi_x(x^{\bar{u}}(t_f))y^{\bar{x},d}(t_f) + c \left [\max_{i\in E} \left |(g^1_i)_x(x^{\bar{u}}(t_f)) y^{\bar{x},d}(t_f)\right |, 
		\right .}
	\nonumber \\
	&{\displaystyle
		\left .
		(g^2_{j_k})_x(x^{\bar{u}}(t_f)) y^{\bar{x},d}(t_f) \right ].}
	\label{5.15}
	\end{eqnarray}
	Extract a subsequence (we do not relabel) such that $j_k= \bar{j}$ for all $k$ and
	$t_k\rightarrow \bar{t}$ for some index value $\bar{j}$.
	By continuity
	of the functions involved $\bar{j}\in I_{0,\bar{u}}$ 
	and (\ref{5.15}) is valid with $\bar{j}$. 
	
	We have arrived at
	\begin{eqnarray}
	&{\displaystyle
		0 \leq  \phi_x(x^{\bar{u}}(t_f)) y^{\bar{x},d}(t_f)
		+c \max \left [ \max_{i\in E}
		\left |(g^1_i)_x(x^{\bar{u}}(t_f))y^{\bar{x},d}(t_f)\right |, \right .}
	\nonumber \\
	&{\displaystyle
		\left .
		\max_{j\in I_{0,\bar{u}}}(g^2_j)_x(x^{\bar{u}}(t_f)) y^{\bar{x},d}(t_f) \right ] .}\label{AppendixFinalInequality} 
	\end{eqnarray}
	\label{nondescfun2}
	Since $\sigma_{c_k}(u_k)\rightarrow 0$, from (\ref{SigmaPos}) we must also have $d_k\rightarrow 0$. $\{c_k\}$ is bounded and $M(u_k)\rightarrow 0$ then (\ref{SigmaPos}) also implies that $\beta_k\rightarrow 0$.
	Therefore, the inequality (\ref{AppendixFinalInequality}), which holds for all $d\in {\mathcal D}_{\bar{u}}$, 
	in particular for $(d\equiv 0)\in {\mathcal D}_{\bar{u}}$, is what we set out to prove.
	
	Finally, we must attend to the case when {\it Algorithm} generates a finite sequence
	which terminates at a control $u_{\bar{k}} = \bar{u}$ satisfying the stopping criterion.
	This case is dealt with by applying the preceding arguments to the infinite
	sequence of controls obtained by `filling in' with repetitions of the following 	control $u_k$. 
\end{proof}    

\section{Derivation of adjoint equations for Case 3--1}
\label{SecAppendixC}

{\it Case 3--1).} In the second case we assume that in the time interval $ [t_0,t_t] $ the system evolves according to DAEs  (17)-(18). At a transition time $ t_t $ the continuous state trajectory leaves the switching surface with the condition $ \eta(x(t_t),u(t_t)) = \alpha(x(t_t),u(t_t)) = 0 $ satisfied---see (10). After the transition the system evolves according to the equation $ x' = f_1(x,u) $ up to an ending time $ t_f $.

To derive the adjoint equations we construct the following augmented functional
\begin{eqnarray}
&{\displaystyle \Phi(x,z,u,\lambda_f,\lambda_h,\pi) = \phi(x(t_f)) + \pi \eta(x(t_t^-),u(t_t^-)) + }  \nonumber\\
&{\displaystyle \int_{t_0}^{t_t^-} \left [\lambda_f^T(t) \left( x'(t) - f_F(x(t),u(t)) - h_x^T(x(t))z(t) \right) + \lambda_h^T(t) h(x(t))\right ] dt}\nonumber \\
&{\displaystyle +\int_{t_t^+}^{t_f} \lambda_f^T(t) \left( x'(t) - f_1(x(t),u(t)) \right) dt.}  \nonumber
\end{eqnarray}
We now calculate the variation of the augmented functional 
\begin{eqnarray}
&{\displaystyle d\Phi(x,z,u,\lambda_f,\lambda_h,\pi) = \phi_x(x(t_f))dx(t_f) + \pi \eta_x(x(t_t^-),u(t_t^-))dx(t_t) + }\nonumber \\
&{\displaystyle \pi \eta_u(x(t_t^-),u(t_t^-))du(t_t^-) +  \lambda_f^T(t_t^-) \left( x'(t_t^-) - f_F(x(t_t^-),u(t_t^-)) - \right .}\nonumber \\
&{\displaystyle \left . h_x^T(x(t_t^-))z(t_t^-) \right)dt_t  +\lambda_h^T(t_t^-) h(x(t_t^-)) dt_t + }\nonumber \\
&{\displaystyle d \left [\int_{t_0}^{t_t^-} \left [\lambda_f^T(t) \left( x'(t) - f_F(x(t),u(t)) - h_x^T(x(t))z(t) \right) + \lambda_h^T(t) h(x(t))\right ]dt \right ] }\nonumber \\
&{\displaystyle 
	-\lambda_f^T(t_t^+) \left( x'(t_t^+) -f_1(x(t_t^+),u(t_t^+)) \right) dt_t}  \nonumber \\
&{\displaystyle + d \left [\int_{t_t^+}^{t_f} \lambda_f^T(t) \left( x'(t) - f_1(x(t),u(t)) \right )dt\right ].} \nonumber
\end{eqnarray}
By taking into account the fact that $dx(t_f) = y^{x,d}(t_f)$
and by integrating by parts the formulas $ \int \lambda_f(t)x(t)dt $ we obtain 
\begin{eqnarray}
&{\displaystyle d\Phi(x,z,u,\lambda_f,\lambda_h,\pi) = \phi_x(x(t_f))y^{x,d}(t_f) +}\nonumber \\
&{\displaystyle  \pi \eta_x(x(t_t^-),u(t_t^-))dx(t_t) + \pi \eta_u(x(t_t^-),u(t_t^-))du(t_t^-) + }\nonumber \\
&{\displaystyle \lambda_f^T(t_t^-) \left( x'(t_t^-) - f_F(x(t_t^-),u(t_t^-)) - h_x^T(x(t_t^-))z(t_t^-) \right)dt_t }\nonumber\\
&{\displaystyle + \lambda_h^T(t_t^-) h(x(t_t^-)) dt_t + d \left [ \left[\lambda_f^T(t)x(t)\right]_{t_0}^{t_t^-}\right ]-}\nonumber \\
&{\displaystyle  d \left [ \int_{t_0}^{t_t^-}\left ( (\lambda_f^T)'(t)  x(t)  + \lambda_f^T(t) \left( f_F(x(t),u(t))  \right . \right . \right .}\nonumber \\
&{\displaystyle \left . \left . \left .+ h_x^T(x(t))z(t) \right)  - \lambda_h^T(t) h(x(t))\right )dt \right ]} \nonumber \\
&{\displaystyle -\lambda_f^T(t_t^+) \left( x'(t_t^+) - f_1(x(t_t^+),u(t_t^+)) \right) dt_t + d \left [ \left[\lambda_f^T(t)x(t)\right]_{t_t^+}^{t_f} \right ]} \nonumber\\
&{\displaystyle -d \left [ \int_{t_t^+}^{t_f} \left ((\lambda_f^T)'(t)  x(t) + \lambda_f^T(t) f_1(x(t),u(t))\right )dt \right ].} \nonumber
\end{eqnarray} 
Expanding further the variations and taking into account the initial conditions of the linearized equations we obtain
\begin{eqnarray}
&{\displaystyle d\Phi(x,z,u,\lambda_f,\lambda_h,\pi) =
	\phi_x(x(t_f))y^{x,d}(t_f) +}\nonumber \\
&{\displaystyle \pi \eta_x(x(t_t^-),u(t_t^-))dx(t_t) +\pi \eta_u(x(t_t^-),u(t_t^-))du(t_t^-) }\nonumber\\
&{\displaystyle +\lambda_f^T(t_t^-)x'(t_t^-)dt_t - \lambda_f^T(t_t^-)f_F(x(t_t^-),u(t_t^-))dt_t} \nonumber \\
&{\displaystyle -\lambda_f^T(t_t^-)h_x^T(x(t_t^-))z(t_t^-)dt_t + \lambda_h^T(t_t^-)h(x(t_t^-))dt_t +}\nonumber \\
&{\displaystyle  \lambda_f^T(t_t^-)y^{x,d}(t_t^- -\int_{t_0}^{t_t^-}\left ( (\lambda_f^T)'(t) y^{x,d}(t) +  \lambda_f^T(t) (f_F)_x(x(t),u(t)) y^{x,d}(t)\right . }  \nonumber\\ 
&{\displaystyle  + \lambda_f^T(t) (f_F)_u(x(t),u(t))d(t) + \lambda_f^T(t) h_x^T(x(t))y^{z,d}(t)  + }\nonumber \\ 
&{\displaystyle \left . \lambda_f^T(t) (h_x^T(x(t))z(t))_x y^{x,d}(t) -\lambda_h^T(t) h_x(x(t))y^{x,d}(t)\right ) dt - \lambda_f^T(t_t^+)x'(t_t^+)dt_t +}\nonumber \\
&{\displaystyle  \lambda_f^T(t_t^+)f_1(x(t_t^+),u(t_t^+))dt_t  +\lambda_f^T(t_f)y^{x,d}(t_f) - \lambda_f^T(t_t^+)y^{x,d}(t_t^+)} \nonumber\\
&{\displaystyle -\int_{t_t^+}^{t_f}\left ( (\lambda_f^T)'(t) y^{x,d}(t) + \lambda_f^T(t) (f_1)_x(x(t),u(t))y^{x,d}(t) +\right .}\nonumber \\
&{\displaystyle \left . \lambda_f^T(t) (f_1)_u(x(t),u(t)) d(t) \right ) dt.} \nonumber 
\end{eqnarray}  
Now we can utilize the formula for the differential $dx(t_t) $, $du(t_t)$ ($du(t_t^-) = d(t_t^-) + u'(t_t^-)dt_t$)
and rearrange the components with respect to differentials $ dx(t_t)$, $ du(t_t)$, $ dt_t $ and variations $ y^{x,d}(t_f)$,  $y^{x,d}(t)$, $ y^{z,d}(t) $, $ d(t)$ to obtain
\begin{eqnarray}
&{\displaystyle d\Phi(x,z,u,\lambda_f,\lambda_h,\pi) =  
	\left( \phi_x(x(t_f)) +\lambda_f^T(t_f) \right)y^{x,d}(t_f) +} \nonumber \\
&{\displaystyle  \left( \pi \eta_x(x(t_t^-), u(t_t^-)) + \lambda_f^T(t_t^-) -\lambda_f^T(t_t^+) \right) dx(t_t)}\nonumber \\
&{\displaystyle + \left (\pi \eta_u(x(t_t^-),u(t_t^-))u'(t_t^-) +   \lambda_f^T(t_t^+)f_1(x(t_t^+),u(t_t^+)) \right . } \nonumber \\ 
&{\displaystyle \left . - \lambda_f^T(t_t^-)f_F(x(t_t^-),u(t_t^-))  - \lambda_f^T(t_t^-)h_x^T(x(t_t^-))z(t_t^-)\right ) dt_t } \nonumber \\
&{\displaystyle  + \lambda_h^T(t_t^-)h(x(t_t^-))  dt_t - \int_{t_0}^{t_t^-} \left (\left( (\lambda_f^T)'(t) +  \lambda_f^T(t) (f_F)_x(x(t),u(t))\right . \right . }\nonumber\\ 
&{\displaystyle + \left. \lambda_f^T(t) (h_x^T(x(t))z(t))_x -\lambda_h^T(t) h_x(x(t)) \right) y^{x,d}(t)}\nonumber \\
&{\displaystyle \left . + \lambda_f^T(t) (f_F)_u(x(t),u(t)) d(t) + \lambda_f^T(t) h_x^T(x(t))y^{z,d}(t)\right ) dt } \nonumber \\
&{\displaystyle -\int_{t_t^+}^{t_f} \left (\left( (\lambda_f^T)'(t) + \lambda_f^T(t) (f_1)_x(x(t),u(t)) \right) y^{x,d}(t)\right . }\nonumber \\
&{\displaystyle \left . + \lambda_f^T(t) (f_1)_u(x(t),u(t))\right ) d(t) dt + \pi \eta_u(x(t_t^-),u(t_t^-))d(t_t^-).} \nonumber 
\end{eqnarray} 

We can now state conditions for adjoint equations such that the expressions with differentials $ dx(t_t),\ dt_t $ and variations $ y^{x,d}(t_f),\ y^{x,d}(t),\ y^{z,d}(t) $ disappear, so eventually only the coefficients with variation $ d(t) $ remain. 

To this end we want to zero the following components
\begin{eqnarray}
&{\displaystyle \left( (\lambda_f^T)'(t) + \lambda_f^T(t) (f_F)_x(x(t),u(t)) + \right .}\nonumber \\
&{\displaystyle \left . \lambda_f^T(t) (h_x^T(x(t))z(t))_x -\lambda_h^T(t) h_x(x(t)) \right) y^{x,d}(t),\ t\in[t_0,t_t^-]}\nonumber \\
&{\displaystyle \lambda_f^T(t) h_x^T(x(t))y^{z,d}(t),\ t\in[t_0,t_t^-] }\nonumber \\
&{\displaystyle \left( (\lambda_f^T)'(t) + \lambda_f^T(t) (f_1)_x(x(t),u(t)) \right) y^{x,d}(t),\ t\in[t_t^+,t_f]. }\nonumber
\end{eqnarray}
We achieve that by assuming that
\begin{eqnarray}
(\lambda_f^T)'(t) & = & -\lambda_f^T(t) (f_F)_x(x(t),u(t)) - \lambda_f^T(t) (h_x^T(x(t))z(t))_x \nonumber \\
& &  + \lambda_h^T(t) h_x(x(t)) \label{eqAdjEqsFirstCaseAdjDaeDiffa}\\ 
0 & = & \lambda_f^T(t) h_x^T(x(t)),\ t\in[t_0,t_t^-) \label{eqAdjEqsFirstCaseAdjDaeAlga}
\end{eqnarray}
and 
\begin{equation}
(\lambda_f^T)'(t) = - \lambda_f^T(t) (f_1)_x(x(t),u(t)),\ t\in[t_t^+,t_f).  
\label{eqAdjEqsFirstCaseAdjOdea}
\end{equation}
with the terminal condition
\begin{equation}
\lambda_f(t_f) = -\phi_x^T(x(t_f)).\nonumber
\end{equation}

To calculate the consistent values of $ \lambda_f(t_t^-) $ and $ \lambda_h(t_t^-) $ the following system of equations have to be solved for the variables $ \lambda_f^-,\ \lambda_h^-,\ \pi,\ \nu_t $ 
\begin{eqnarray}
&{\displaystyle \lambda_f^- - \lambda_f(t_t^+) + \pi \eta_x^T(x(t_t^-),u(t_t^-)) =  \nu_t h_x^T(x(t_t^-)) } \label{CIC31a} \\
&{\displaystyle \left (\lambda_f^-\right )^T f_F(x(t_t^-),u(t_t^-))  +\left (\lambda_f^-\right )^Th_x^T(x(t_t^-))z(t_t^-)  }\nonumber \\
&{\displaystyle - \lambda_h^-h(x(t_t^-)) - \pi \eta_u(x(t_t^-),u(t_t^-))u'(t_t^-) = \lambda_f^T(t_t^+) f_1(x(t_t^+),u(t_t^+))} \label{CIC31b} 
\end{eqnarray}
and
\begin{eqnarray}
0 & = & h_x(x(t_t^-))\lambda_f^-\label{CIC31c}  \\
0 & = &\left(h_x(x(t_t^-))\right)'\lambda_f^- -  \nonumber \\
& &  h_x(x(t_t^-)) \left( f_F \right)_x^T(x(t_t^-),u(t_t^-))\lambda_f^- - \nonumber \\
& & h_x(x(t_t^-)) \left(h_x^T(x(t_t^-))z(t_t^-)\right)_x^T\lambda_f^- +  \nonumber \\
& &  h_x(x(t_t^-)) h_x^T(x(t_t^-)) \lambda_h^- \label{CIC31d} 
\end{eqnarray}
where $ \nu_t $ is a number.

One can solve these equations in the following way. First of all, we multiply the equation (\ref{CIC31a}) by $f_F(x(t_t^-),u(t_t^-)) + h_x(x(t_t^-))^T z(t_t^-)$ and take into account (\ref{CIC31b}), that $h_x(x(t_t^-))\left (f_F(x(t_t^-),u(t_t^-)) + h_x(x(t_t^-))^T z(t_t^-)\right ) = 0$, $h(x(t_t^-))=0$, to obtain (which is possible due to ${\bf (H2)}$)
\begin{eqnarray}
&{\displaystyle \pi = \lambda_f(t_t^+)^T\left (f_F(x(t_t^-),u(t_t^-)) + h_x(x(t_t^-))^T z(t_t^-) - f_1(x(t_t^+),u(t_t^+))\right )/}\nonumber \\
&{\displaystyle \left [\eta_x(x(t_t^-),u(t_t^-))\left ( f_F(x(t_t^-),u(t_t^-)) + h_x(x(t_t^-))^T z(t_t^-)\right ) + \right . }\nonumber \\
&{\displaystyle \left . \eta_u(x(t_t^-),u(t_t^-))u'(t_t^-) \right ]} \nonumber
\end{eqnarray}
Next, we multiply (\ref{CIC31a}) by $h_x(x(t_t^-))$ and take into account (\ref{CIC31c}) to obtain (which is possible due to ${\bf (H1)}$)
\begin{eqnarray}
&{\displaystyle \nu_t = h_x(x(t_t^-)) \left (- \lambda_f(t_t^+) + \pi \eta_x(x(t_t^-),u(t_t^-))^T    \right ) /\left \| h_x(x(t_t^-))\right \|^2}\nonumber
\end{eqnarray}
Having $\pi$ and $\nu_t$ we can evaluate then 
\begin{eqnarray}
{\displaystyle   \lambda_f^-  = \lambda_f(t_t^+) - \pi \eta_x^T(x(t_t^-),u(t_t^-)) +  \nu_t h_x^T(x(t_t^-)) }\nonumber 
\end{eqnarray}
and eventually
\begin{eqnarray}
&{\displaystyle  \lambda_h^-  = \left (\lambda_f^-\right )^T \left ( \left ( \left( f_F \right)_x(x(t_t^-),u(t_t^-)) +  \left(h_x(x(t_t^-))z(t_t^-)\right)_x\right ) h_x^T(x(t_t^-)) - \right . }\nonumber \\
&{\displaystyle \left .  \left(h_x^T(x(t_t^-))\right)'  \right )/ \left \| h_x(x(t_t^-))\right \|^2.}\nonumber 
\end{eqnarray}
Next we set $\lambda_f(t_t^-) = \lambda_f^-$, $\lambda_h(t_t^-)=\lambda_h^-$.

Once we solve the adjoint equations we obtain the adjoint variables $ \lambda_f $ and $ \lambda_h $. Now we can calculate the first variation of a cost function $ \phi(t_f) $ with respect to a control function variation $ d $ as follows
\begin{eqnarray}
&{\displaystyle d\phi(x(t_f)) = d\Phi(x,z,u,\lambda_f,\lambda_h,\pi) = \pi \eta_u(x(t_t),u(t_t^-)) d(t_t^-)} \nonumber\\
&{\displaystyle -\int_{t_0}^{t_t^-} \lambda_f^T(t) (f_F)_u(x(t),u(t)) d(t)  dt - \int_{t_t^+}^{t_f} \lambda_f^T(t) (f_1)_u(x(t),u(t)) d(t) dt.} \nonumber 
\end{eqnarray} 

\section{The weak maximum principle for {\it Case 3--1}}
\label{SecAppendixD}

\noindent ${\bf (NC^{31})}$: There exist: nonnegative numbers $\alpha_j^2$, $j\in I$, numbers $\alpha_i^1$, $i\in E$ such that $\sum_{i\in E}\left | \alpha_i^1\right | +\sum_{j\in I} \alpha_j^2 \neq 0$; piecewise smooth function $\lambda_f$; piecewise continuous function $\lambda_h$, such that the following hold\\
{\it (i) terminal conditions}
\begin{eqnarray}
&{\displaystyle 
-\lambda_f(t_f) = \phi_x^T(\bar{x}(t_f)) + \sum_{i\in E} \alpha_i^1\left (g_i^1\right )_x^T (\bar{x}(t_f)) +  \sum_{j\in I} \alpha_j^2\left (g_j^2\right )_x^T (\bar{x}(t_f))}\nonumber
\end{eqnarray}
{\it (ii) adjoint equations}\\
for $ t\in [t_t,t_f) $
\begin{eqnarray}
\lambda_f' = - \left( f_1 \right)_x^T(\bar{x},\bar{u})\lambda_f;
\nonumber 
\end{eqnarray}
\noindent for $ t\in [t_0,t_t) $
\begin{eqnarray}
\lambda_f' &=& - \left( f_F \right)_x^T(\bar{x},\bar{u})\lambda_f - \left(h_x^T(\bar{x})\bar{z}\right)_x^T\lambda_f + h_x^T(\bar{x}) \lambda_h \nonumber \\ 
0 &=& h_x(\bar{x})\lambda_f \nonumber
\end{eqnarray}
{\it (iii) jump conditions}---$\lambda_f(t_t^-)$, $\lambda_h(t_t^-)$, $\nu_t$ and $\pi$ are evaluated from the equations below
\begin{eqnarray}
&{\displaystyle 
\lambda_f(t_t^-) - \lambda_f(t_t^+) + \pi \eta_x^T(\bar{x}(t_t),\bar{u}(t_t))  =  \nu_t h_x^T(\bar{x}(t_t))  }\nonumber \\
&{\displaystyle \lambda_f^T(t_t^-) f_F(\bar{x}(t_t),\bar{u}(t_t)) + \lambda_f^T(t_t^-)h_x^T(\bar{x}(t_t))\bar{z}(t_t^-) } \nonumber \\
&{\displaystyle - \pi \eta_u(\bar{x}(t_t),\bar{u}(t_t))\bar{u}'(t_t^-) - \lambda_h^T(t_t^-)h(\bar{x}(t_t))  =}\nonumber \\
&{\displaystyle \lambda_f^T(t_t^+) f_1(\bar{x}(t_t),\bar{u}(t_t^+)) } \nonumber
\end{eqnarray}
and
\begin{eqnarray}
0 & = & h_x(\bar{x}(t_t))\lambda_f(t_t^-) \nonumber \\
0 & = &\left(h_x(\bar{x}(t_t))\right)'\lambda_f(t_t^-) - h_x(\bar{x}(t_t)) \left( f_F \right)_x^T(\bar{x}(t_t),\bar{u}(t_t))\lambda_f(t_t^-) - \nonumber \\
& & h_x(\bar{x}(t_t)) \left(h_x^T(\bar{x}(t_t))\bar{z}(t_t^-)\right)_x^T\lambda_f(t_t^-) +  h_x(\bar{x}(t_t)) h_x^T(\bar{x}(t_t)) \lambda_h (t_t^-). \nonumber
\end{eqnarray}
{\it (iv) the weak maximum principle}\\
\begin{eqnarray}
    &{\displaystyle H^{31}(\bar{x},\bar{u},\lambda_f,t_t,\pi ,u)\leq H^{31}(\bar{x},\bar{u},\lambda_f,t_t, \pi, \bar{u}),\ \forall u\in {\mathcal U}}\nonumber 
\end{eqnarray}
where $H^{31}(\bar{x},\bar{u},\lambda_f,t_t, \pi, u) = \int_{t_0}^{t_t^-}\lambda_f^T(t) \left (f_F\right )_u (\bar{x}(t),\bar{u}(t))u(t)dt +$\\ $\int_{t_t^+}^{t_f}\lambda_f^T(t) \left (f_1\right )_u (\bar{x}(t),\bar{u}(t))u(t)dt -\pi \eta_u (\bar{x}(t_t),\bar{u}(t_t))u(t_t)$.\\
{\it (v) complementarity conditions}
\begin{eqnarray}
\alpha_j^2 = 0,\ \ {\rm if}\ \ j\not\in I_{0,\bar{u}}.
\nonumber 
\end{eqnarray}
{\it (vi) feasibility} 
\begin{eqnarray}
&{\displaystyle 
g_i^1 (\bar{x}(t_f)) = 0,\  i\in E,\  g_j^2 (\bar{x}(t_f))\leq 0,\ j\in I.}\nonumber 
\end{eqnarray}